\documentclass[leqno,12pt]{amsart} 
\usepackage{amssymb,amscd}
\theoremstyle{plain} 
\newtheorem{theorem}{\indent\sc Theorem}[section]
\newtheorem{lemma}[theorem]{\indent\sc Lemma}
\newtheorem{corollary}[theorem]{\indent\sc Corollary}
\newtheorem{proposition}[theorem]{\indent\sc Proposition}

\newtheorem*{maintheorem}{\indent\sc Theorem}
\theoremstyle{definition} 

\newtheorem{remark}[theorem]{\indent\sc Remark}
\newtheorem{example}[theorem]{\indent\sc Example}

\newcommand\mylabel[1]{\label{#1}}

\newcommand{\ZZ}{\mathbb{Z}}
\newcommand{\QQ}{\mathbb{Q}}

\newcommand{\FF}{\mathbb{F}}

\newcommand{\PP}{\mathbb{P}}

\newcommand{\GG}{\mathbb{G}}

\newcommand  {\shExt}   {\mathcal{E} \!\text{\textit{xt}}}

\newcommand  {\shF}     {\mathcal{F}}

\newcommand  {\shHom}   {\mathcal{H}\!\text{\textit{om}}}
\newcommand  {\shI}     {\mathcal{I}}

\newcommand  {\shK}     {\mathcal{K}}
\newcommand  {\shM}     {\mathcal{M}}

\newcommand  {\shN}     {\mathcal{N}}
\newcommand  {\shL}     {\mathcal{L}}

\newcommand  {\shT}     {\mathcal{T}}

\newcommand  {\shProj}  {\mathcal{P}\!\text{\textit{roj}}\,}

\newcommand  {\Aut}     {\operatorname{Aut}}

\newcommand  {\cid}     {\mathfrak{c}}

\renewcommand{\cong}    {\equiv}

\newcommand  {\cris}    {{\operatorname{cris}}}

\newcommand  {\Hom}     {\operatorname{Hom}}

\newcommand  {\id}      {\operatorname{id}}

\newcommand  {\Isom}    {\operatorname{Isom}}

\newcommand  {\lra}     {\longrightarrow}

\newcommand  {\maxid}   {\mathfrak{m}}

\newcommand  {\mt}      {\mapsto}

\renewcommand{\O}       {\mathcal{O}}

\newcommand  {\Pic}     {\operatorname{Pic}}

\newcommand  {\pr}      {\operatorname{pr}}

\newcommand  {\quadand} {\quad\text{and}\quad}

\newcommand  {\ra}      {\rightarrow}

\newcommand  {\Sing}    {\operatorname{Sing}}

\newcommand  {\Spec}    {\operatorname{Spec}}

\begin{document}

\title[Weak del Pezzo surfaces]{Weak del Pezzo surfaces with   irregularity} 

\author[S.\ Schr\"oer]{Stefan Schr\"oer} 



%
\address{
Mathematisches Institut\endgraf
Heinrich-Heine-Universit\"at\endgraf
40225 D\"usseldorf\endgraf
Germany
}
\email{schroeer@math.uni-duesseldorf.de}


\maketitle

\begin{abstract}
I construct   normal 
del Pezzo surfaces, and regular weak del Pezzo surfaces
as well, with positive irregularity $q>0$.
This can happen only over nonperfect fields.
The surfaces in question are twisted forms of nonnormal del Pezzo surfaces,
which were classified by Reid.
The twisting is with respect to the flat topology and infinitesimal
group scheme actions. The twisted surfaces appear as generic
fibers for   Fano-Mori contractions on certain   threefolds
with only canonical singularities.
\end{abstract}

\section*{Introduction} 
\mylabel{Introduction}

Suppose that $X$ is a smooth and projective scheme over the complex numbers.
The   Kawamata-Viehweg Vanishing Theorem 
asserts that $H^i(X,\omega_X\otimes\shL)=0$ for all  integers $i>0$
and all invertible sheaves $\shL$ that are nef and big (see \cite{Kawamata 1982} and \cite{Viehweg 1982}).
In the special case that $X$ is a \emph{weak Fano variety}, in other words, the dual of the dualizing sheaf
is nef and big, we may apply the Kawamata-Viehweg Vanishing Theorem with $\shL=\omega_X^\vee$
and conclude that  $H^i(X,\O_X)=0$ for all integers $i>0$. 
\emph{It is unknown whether or to what extend this particular vanishing holds true
for Fano or weak Fano varieties in positive characteristics.} 

The Kawamata-Viehweg Vanishing Theorem is a generalization
of the Kodaira Vanishing Theorem \cite{Kodaira  1953}, which deals with
ample rather than nef and big invertible sheaves.
It is well-known that the Kodaira Vanishing Theorem
does not hold true in positive characteristics. Raynaud 
\cite{Raynaud 1978} constructed the first counterexamples,
which are fibered surfaces whose generic fiber is regular but not smooth. The surfaces are mostly of general type.
 A rather different set of counterexamples is due to 
Lauritzen \cite{Lauritzen 1996} relying on representation theory: He used homogeneous schemes of the form $G/B$, where $B\subset G$ is a nonreduced
Borel subgroup scheme in some  linear algebraic group. Using more elementary methods,
Lauritzen and Rao \cite{Lauritzen; Rao 1997}  further  constructed smooth Fano varieties of dimension $d\geq 6$
so that Kodaira vanishing fails for some ample invertible sheaves $\shL\neq\omega_X^\vee$.

Esnault  \cite{Esnault 2003} gives a completely different aspect involving crystalline cohomology:
Her results, which apply to a much wider class than just Fano varieties,  tell
us that for smooth Fano varieties over perfect fields $k$,
with the ring of Witt vectors $W$ and field of fractions $W\subset K$,
the following holds: The   part with slopes $\lambda\in [0,1[$ inside the
crystalline cohomology groups
$H_\cris ^i(X/W)\otimes_W K$ vanishes for $i>0$. 
On the other hand, this part $(H_\cris ^i(X/W)\otimes_W K)_{[0,1[}$  is isomorphic to Serre's Witt vector cohomology $H^i(X,W\O_X)\otimes_W K$. In turn, the group $H^i(X,W\O_X)$    is related to ordinary cohomology
groups $H^i(X,\O_X)$ by exact sequences, but it seems difficult to gain control over
torsion phenomena. Note that Berthelot, Bloch and Esnault
\cite{Berthelot; Bloch; Esnault 2005} extended the bijection between the slope $[0,1[$-part and Witt vector
cohomology to singular schemes, with rigid cohomology instead of crystalline 
cohomology.

There are some positive results in low dimensions. If follows
from the classification of smooth   del Pezzo surfaces, which are the 2-dimensional   Fano varieties,
that $H^1(S,\O_S)=0$ holds regardless to the characteristic. The same
holds for weak del Pezzo surfaces. For a nice account, see \cite{Demazure 1976}.
Shepherd-Barron  \cite{Shepherd-Barron 1997} established the vanishing $H^1(X,\O_X)=H^2(X,\O_X)=0$
for smooth Fano threefolds.
On the other hand, Reid \cite{Reid 1994} constructed nonnormal del Pezzo surfaces
with $H^1(S,\O_S)\neq 0$.

My original motivation for this work was to construct  \emph{regular} Fano varieties
over  nonperfect fields that   have  $H^1(X,\O_X)\neq 0$. 
The point here is that regularity does not imply geometric regularity (= formal smoothness) over
nonperfect fields.
I did not quite succeed in my goals, but I came close to it.
The main result of this paper is as follows:

\begin{maintheorem}
Over every nonperfect field of  characteristic $p=2$, there are weak del Pezzo surfaces 
that are regular, and normal  del Pezzo surfaces $S$ with  only factorial
rational double points of type $A_1$ as formal singularities, both with $h^1(\O_S)\neq 0$.
\end{maintheorem}

These del Pezzo surfaces are   twisted forms of  Reid's nonnormal del Pezzo
surfaces, and the weak del Pezzo surface is obtained by resolving the singularity.
Such   del Pezzo surfaces  necessarily become nonnormal after passing to the perfect closure of the ground field.
Indeed, it follows from the work of Hidaka and Watanabe \cite{Hidaka; Watanabe 1981} and myself \cite{Schroeer 2001}
that $h^1(\O_S)=0$ for   geometrically normal del Pezzo surfaces.

The existence of such wild del Pezzo surfaces $S$ over nonperfect fields   has   consequences
for the structure theory of algebraic varieties $X$ over algebraically closed fields.
Namely, such del Pezzo surfaces might arise as generic fibers in some  Fano-Mori contractions
of fiber type, obtained by contracting an extremal ray.  
To my knowledge, the geometry of fibrations $f:X\ra B$ whose generic fiber $X_\eta$
is   not geometrically regular or geometrically normal
has not been studied systematically, except for the  quasielliptic surfaces
of Bombieri and Mumford, see \cite{Bombieri; Mumford 1977} and \cite{Bombieri; Mumford 1976}.

The existence of Fano-Mori contractions of extremal rays on smooth threefolds
in arbitrary characteristics was established by Koll\'ar in \cite{Kollar 1991}.
In Remark 1.2, he raised the question whether there are contractions of fiber type
whose generic geometric fibers are  nonnormal del Pezzo surfaces.
We shall see that our exotic del Pezzo surfaces appear as generic fibers $S=Z_\eta$
for some Fano-Mori contraction $f:Z\ra E$ of fiber type, where $Z$ is a threefold,
and $E$ is the  supersingular elliptic curve in characteristic two.
Unfortunately, my results are not strong enough to make the total space smooth.
However, the threefold $Z$ will be locally of complete intersection, locally factorial, 
with only canonical singularities. The anticanonical divisor is nef and has Kodaira dimension two, 
and the first higher direct image  $R^1f_*(\O_Z)$ 
is nonzero.

There are several papers dealing with Fano threefolds in positive characteristics.
For example, Shepherd-Barron \cite{Shepherd-Barron 1997} 
obtained a classification for Picard number $\rho=1$,
and Megyesi \cite{Megyesi 1998} treated the case of Fano varieties of index $\geq 2$.
Saito showed that on Fano threefolds with Picard number $\rho=2$
there are no fibrations whose geometric generic fiber is a nonnormal del Pezzo surface
\cite{Saito 2003}. Mori and Saito \cite{Mori;  Saito 2003} have further results
 on wild hypersurface bundles.

Here is a plan for the paper:
In Section \ref{Twisted forms}, we collect some general facts
on twisting and twisted forms, and give a criterion for
regularity of twisted forms.
In Section \ref{Ribbons}, we recall Reid's construction of  
nonnormal del Pezzo surfaces $Y$ in terms of glueing
along a double line  to a rational cuspidal curve, and 
discuss the glueing process in detail.
In Section \ref{The geometric construction},
we shall see that the resulting del Pezzo surface $Y$ is
locally of complete intersection.
In Section \ref{Picard scheme and dualizing sheaf}, I analyse
the Picard group and the dualizing sheaf on these del Pezzo surfaces.
Section \ref{Cartier divisors and Weil divisors} contains
a discussion of curves of degree one. In particular, we shall
explain how and why some of these curves are Cartier divisors, and others are
only Weil divisors.
In Section \ref{The tangent sheaf}, we shall prove that
the cotangent sheaf modulo torsion is locally free of rank two.
This  seems to be a rather
special situation. An immediate consequence is that the
tangent sheaf is locally free. Under suitable assumptions,
we moreover identify a global vector field $\delta\in H^0(Y,\Theta_{Y/k})$
that defines an $\alpha_2$-action.
Unfortunately, this vector field has a unique zero at the so-called
point at infinity $y_\infty\in Y$.
However, we check in Section \ref{Splitting type of tangent sheaf}
that zeros of vector fields are inevitable, by computing the
splitting type of the tangent sheaf. It turns out that our choice of $\delta$
is in some sense the best possible.
In Section \ref{Twisted del Pezzo surfaces}, we use the
$\alpha_2$-action to construct twisted forms $Y'$ of our
del Pezzo surface $Y$. It turns out that the twisted forms
are normal, with only one singularity. A formal analysis
reveals that this singularity is a rational double point of
type $A_1$, which is moreover factorial.
In Section \ref{Fano--Mori contractions of fiber type} we use our results
to construct some interesting Fano-Mori contractions.
In Section \ref{Maps to projective spaces}, we finally discuss
ampleness and semiampleness for line bundles on $Y$.
Here we see that it is impossible to realize $Y$ as
 hypersurface  or   double covering of  some $\PP^n$.
The smallest embeddings have codimension three, and the smallest coverings have degree four.

\medskip
I wish to thank Torsten Ekedahl, H\'el\`ene Esnault, and Christian Liedtke for stimulating discussions, 
and Burt Totaro for helpful comments.

\section{Twisted forms}
\mylabel{Twisted forms}

In this section I discuss some useful general aspects of twisting and twisted form
that we shall apply later to nonnormal del Pezzo surfaces.
Suppose $S$ is a scheme of finite type over a base field $k$.
Another $k$-scheme $S'$ is called a \emph{twisted form} of $S$
if there is a nonzero $k$-algebra $R$ with $S_R\simeq S'_R$.
Such   schemes $S'$ are automatically of finite type by descent theory
(see \cite{SGA 1}, Expos\'e VIII, Proposition 3.3).

\begin{lemma}
\mylabel{finite extension}
If $S'$ is a twisted form of $S$, then there is a finite
field extension $k\subset E$ with $S_E\simeq S'_E$.
\end{lemma}

\proof
Choose an isomorphism $f:S_R\ra S'_R$.
As explained in \cite{EGA IVc}, Theorem 8.8.2, there is a  $k$-subalgebra of finite
type $R_\alpha\subset R$ and an isomorphism $f_\alpha:S_{R_\alpha}\ra S'_{R_\alpha}$ inducing $f$.
Choose a maximal ideal $\maxid\subset R_\alpha$. Then $E=R_\alpha/\maxid$ is
a finite field extension of $k$, and we may restrict $f_\alpha$ to $E$.
\qed

\medskip
If follows that the set of isomorphism classes of twisted forms $S'$ is a subset of 
the nonabelian cohomology set $H^1(k,\Aut_{S/k})$,
where we may use the finite flat topology.
A nice account of this correspondence in the context of Galois cohomology
appears in Serre's book \cite{Serre 1972}. The full theory is exposed at length
in Giraud's treatise \cite{Giraud 1971}.
The basic construction goes as follows: Let $T$ be a torsor under
$\Aut_{S/k}$ with action from the left. Then we may form 
the product $S\times T$ and obtain the  quotient by the diagonal action
$$
S\wedge T=\Aut_{S/k}\backslash(S\times T),\quad (s,t)\sim (gs,gt).
$$
Note that $\Aut_{S/k}$ also acts from the right on $S$ via
$sg=g^{-1}s$, so we may rewrite the equivalence relation in the
particularly attractive form  $(s,t)\sim (sg^{-1},gt)$.

The result $S\wedge T$, which is a sheaf in the finite flat topology,  is a sheaf
twisted form of $S$. This sheaf, however, is not necessarily representable by a scheme.
We shall discuss this below.
Conversely, if $S'$ is a twisted form of $S$, then
$T=\underline{\Isom}(S',S)$ is an $\Aut_{S/k}$-torsor with action from the left, and
the canonical map $S\wedge T\ra S'$, $(s,t)\mapsto t^{-1}(s)$  is an isomorphism.

Now let $G\subset\Aut_{S/k}$ be a subgroup scheme.
Then we have an induced map on nonabelian cohomology $H^1(k,G)\ra H^1(k,\Aut_{S/k})$.
Given any $G$-torsor $T$, we may form $S\wedge T=G\backslash(S\times T)$ to produce twisted forms.
Note that the twisted form  might be trivial,
although the torsor  is nontrivial. More precisely:

\begin{lemma}
\mylabel{trivial twisting}
The twisted form $S\wedge T$ is isomorphic to $S$ if and only if there
is a $G$-equivariant morphism $T\ra\Aut_{S/k}$.
\end{lemma}

\proof
As explained in \cite{Giraud 1971} Chapter III, Proposition 3.2.2 , we have a sequence
$$
H^0(k,\Aut_{S/k})\lra H^0(k,G\backslash\Aut_{S/k})\lra H^1(k,G)\lra H^1(k,\Aut_{S/k}),
$$
which is exact in the following sense: The $G$-torsors $T$ inducing trivial $\Aut_{S/k}$-torsors
come from the sections $x\in G\backslash\Aut_{S/k}$. The $G$-torsor coming from such a section
$x$ is  
the fiber over $x$ under the projection $\Aut_{S/k}\ra G\backslash\Aut_{S/k}$, whence the assertion.
\qed

\medskip
We are mainly interested in the case that the group scheme $G$ is finite.
Then there are almost no problems with representability:

\begin{lemma}
\mylabel{algebraic space}
Suppose $G$ is a finite group scheme. Then $S'=S\wedge T$ is
an algebraic space. It is a scheme if either $G$ is infinitesimal,
or if $S$ is quasiprojective.
\end{lemma}

\proof
It follows from \cite{SGA 1}, Expos\'e VIII, Corollary 7.7,
that the torsor $T$ is representable by a scheme. The quotient
$S'=G\backslash(S\times T)$
exists as an algebraic space, according to very general results 
of Keel and Mori \cite{Keel; Mori 1997}.
If $G$ is infinitesimal or if $S$ is quasiprojective then 
the $G$-invariant affine open subset $U_\alpha\subset S\times T$ form
a covering. By \cite{SGA 3a}, Expos\'e V, Theorem 4.1 the quotient $S'$ of $S\times T$ by the free
$G$-action exists as a scheme.
\qed

\medskip
From now on we assume for convenience that the group scheme $G$ is finite, and that $S'=S\wedge T$ is a scheme,
where $T$ is a $G$-torsor. 

\begin{lemma}
\mylabel{smooth complete}
The scheme $S$ is locally of complete intersection or smooth
if and only if $S'$ is locally of complete interesection or smooth, respectively.
\end{lemma}

\proof
This follows from  \cite{EGA IVd}, Corollary 19.3.4 and Proposition 17.7.1.
\qed

\medskip
In contrast,  regularity or nonregularity does not transfer to twisted forms.
It is possible to remove singularities by passing to twisted forms,
which is  indeed the leitmotiv of this paper. Of course, such   things may happen only
in positive characteristics. I found the following basic fact very useful.

\begin{theorem}
\mylabel{twisted forms}
Let $A\subset S$ be a $G$-invariant
subscheme   whose ideal is locally generated by regular sequences.
If the twisted form $A'=A\wedge T$ is a regular scheme, then the
twisted form $S'=S\wedge T$ is a regular scheme at all points on the subset $A'\subset S'$.
\end{theorem}

\proof
Consider the  commutative diagram
$$
\begin{CD}
A\times T @>>> S\times T\\
@VVV  @VVV\\
A' @>>> S'.
\end{CD}
$$
The vertical maps are surjective and flat, because $G$ acts freely on $A\times T$ and $S\times T$.
Moreover, the diagram is cartesian.
By assumption, the embedding $A\subset S$ is regular.
Hence the induced embedding $A\times T\subset S\times T$ is regular as well.
According to \cite{EGA IVd}, Proposition 19.1.5, this implies that the embedding
$A'\subset S'$ is regular.
By assumption, the scheme $A'$ is regular. If follows that
the scheme $S'$ is regular at all points $s\in A'$, by
\cite{EGA IVd}, Proposition 19.1.1.
\qed

\medskip
We shall mainly apply this in the cases that $A\subset S$ is either a Cartier divisor or
an Artin subscheme. Let me record the latter:

\begin{corollary}
\mylabel{regular ring}
Let $A=Gs$ be the orbit
of a rational point $s\in S$, and let $s'\in S'$ be a closed
point in  $A'\subset S'$. If the scheme $T$ is reduced
and $Gs\subset S$ is a regular embedding,
then the local ring $\O_{S',s'}$ is regular.
\end{corollary}

\proof
The orbit $A=Gs$ of our finite group scheme $G$ is isomorphic to the homogeneous
space $G/H$, where $H=G_s$ is the isotropy group scheme. The projection $G/H\times T\ra H\backslash T$,
$(gH,t)\mapsto Hg^{-1}t$ is well-defined, and induces a bijection
$G\backslash(G/H\times T)\ra H\backslash T$. If follows that  the twisted form
$A'=A\wedge T$ is isomorphic to
$ H\backslash T$. By assumption, the
Artin scheme $T$ is reduced, whence the qotient scheme $H\backslash T$ is reduced as well.
In other words, the twisted form $A'$ is regular, and the Theorem applies.
\qed

\begin{example}
Consider the global field $k=\FF_2(t)$  in characteristic two and the  1-dimensional scheme $S=\Spec k[u^2,u^3]$, which
contains a cuspidal singularity at the origin.
The finite infinitesimal group scheme $G=\alpha_2$ acts on $S$ via the derivation $u^3\mapsto 1$.
It also acts on the   $T= \Spec k[\sqrt{t}]$   via the
derivation $\sqrt{t}\mapsto 1$. The twisted form $S'=S\wedge T$ then must be a regular curve.
Indeed, it is the spectrum of the subalgebra $k[u^2,u^3+\sqrt{t}]\subset k[u^2,u^3,\sqrt{t}]$.
\end{example}

Now suppose that $x\in S$ is a rational point that is fixed under the $G$-action.
Setting $A=\left\{x\right\}$, we see that the twisted form $A'=A\wedge T$ is given by
another rational point $x'\in S'$.  The following tells us that it is impossible
to remove singularities by twisting at points that are  both fixed and singular.

\begin{proposition}
\mylabel{singularities stay}
Assumptions as above. If the local ring $\O_{S,x}$ is not regular, then the local ring
$\O_{S',x'}$ is not regular as well.
\end{proposition}

\proof
Let $B\subset S\times T$ be the preimage of $x\in S$ under the projection map,
which coincides with the preimage of $x'\in S'$ under the quotient map.
Suppose $s'\in S'$ is a regular point. Then the residue field $\kappa(s')$ has
finite projective dimension. Hence $\O_B$ has finite projective dimension as well,
because the quotient map $S\times T\ra S'$ is flat.
Choose  a   resolution 
$
\ldots \ra F_1\lra F_0\ra\kappa(s)\ra 0
$
with finitely generated free $\O_{S,x}$-modules. Pulling back under the flat projection map,
we obtain a free resolution $\ldots \ra F'_1\lra F'_0\ra\O_B\ra 0$.
Whence the kernel of some $F_{i+1}'\ra F_i'$ is free.
By descent theory, the kernel of $F_{i+1}\ra F_i$ must also be free.
In other words, $x\in S$ is regular.
\qed

\section{Glueing along ribbons}
\mylabel{Ribbons}

Throughout the following sections, we shall study the geometry of certain  nonnormal del Pezzo surfaces $Y$ with
irregularity $h^1(\O_Y)>0$.
Such surfaces were first constructed by Reid \cite{Reid 1994}.
A key point in his construction is the use of certain infinitesimal neighborhoods called
ribbons.  Reid's construction works roughly as follows: 

We fix a base field $k$  of characteristic $p=2$.
Let $X=\PP^2$ be the projective plane, and $A=\PP^1$ be the projective line. 
Choose an embedding $A\subset X$ of degree one,
and let $B=A^{(1)}$ be the first order infinitesimal neighborhood, which is a nonreduced quadric.
Let  $C$ be the rational cuspidal curve with arithmetic genus $h^1(\O_C)=1$,
whose normalization is $A\ra C$. The idea now is to 
extend the nonflat normalization map $A\ra C$ to a flat morphism  $\varphi:B\ra C$ of degree two, and 
obtain the desired del Pezzo surface $Y$ via the cocartesian diagram
$$
\begin{CD}
B @>>> X\\
@V\varphi VV @VV\nu V\\
C @>>> \phantom{.}Y.
\end{CD}
$$
Note that the normalization map $\nu:X\ra Y$ is a homeomorphism.
One way to think about this   is that we
thinned out the structure sheaf $\O_X$ by artificially removing     sections
satisfying certain conditions on $B$
to obtain the structure sheaf $\O_Y$,
as explained in Serre's book \cite{Serre 1975}, Chapter IV,  \S1.1.
In some sense, we introduced a curve of cusps $C\subset Y$, which itself contains
a cuspidal singularity.
Naturally, the singularity on this curve of singularities plays a crucial role
in the whole affair.

To make this construction explicit and to explore its properties, it seems inevitable to
introduce coordinates.
Choose indeterminates $u,v$ and cover the projective plane $X=\PP^2$
in the usual way  by three affine open subschemes
\begin{equation}
\label{covering X}
X=\Spec k[u,v]\;\cup\;\Spec k[u^{-1},vu^{-1}]\;\cup\;\Spec k[uv^{-1},v^{-1}].
\end{equation}
We sometimes denote this open affine covering by $X=U\cup U'\cup U''$.
We shall see that our constructions do not work well in homogeneous coordinates, and
it seems   necessary to introduce inhomogeneous coordinates.
The projective line $A=\PP^1$ shall be  
embedded into the projective plane $X=\PP^2$  by setting
\begin{equation}
\label{covering A}
A= \Spec k[u,v]/(v)\;\cup\;\Spec k[u^{-1},vu^{-1}]/(vu^{-1}).
\end{equation}
We write the rational cuspidal curve  $C$  with arithmetic genus $p_a=1$ as the union
of two affine open subschemes
\begin{equation}
\label{covering C}
C=\Spec k[u^2,u^3]\;\cup\;\Spec k[u^{-1}].
\end{equation}
Then we have a canonical morphism $A\ra C$, which is the normalization map.
Finally, consider the first order infinitesimal neighborhood $B=A^{(1)}$ inside 
the projective space $ X=\PP^2$.
This nonreduced quadric  is given by
\begin{equation}
\label{covering B}
B= \Spec k[u,\epsilon]\;\cup\;\Spec k[u^{-1},\epsilon u^{-1}],
\end{equation}
where $\epsilon$ denotes the residue class of $v$ modulo $v^2$.
The inclusion $A\subset B$ is a \emph{ribbon} in the sense
of Bayer and Eisenbud \cite{Bayer; Eisenbud 1995}.
This means that the ideal $\shI\subset\O_B$ of the closed embedding $A\subset B$ satisfies
$\shI^2=0$ and that $\shI$ is an invertible $\O_A$-module.
Note that the first condition implies that the $\O_B$-module structure on $\shI$
indeed comes from an $\O_A$-module structure.
We have an exact cotangent sequence of $\O_A$-modules
\begin{equation}
\label{cotangent sequence}
0\lra\shI\lra\Omega^1_{B/k}\otimes\O_A\lra\Omega^1_{A/k}\lra 0,
\end{equation}
where we use $\shI=\shI/\shI^2$.
Pulling back the extension along the universal derivation $d:\O_A\ra\Omega^1_{A/k}$,
we obtain an extension of sheaves of $k$-vector spaces
$$
0\lra\shI\lra\O_B\lra \O_A\lra 0.
$$
One may recover the multiplication in $\O_B$ by exploiting the fact that
$d:\O_A\ra\Omega^1_{A/k}$ is a derivation. Note that in particular there is
a cartesian diagram of $\O_B$-modules
$$
\begin{CD}
\Omega^1_{B/k}\otimes\O_A@>>>\Omega^1_{A/k}\\
@Ad\otimes 1AA @AAd A\\
\O_B @>>> \phantom{.}\O_A.
\end{CD}
$$
The normalization map $A\ra C$ induces an 
$\O_A$-linear map $\Omega^1_{C/k}\otimes\O_A\ra\Omega^1_{A/k}$.
We may use the latter map to pull back the extension (\ref{cotangent sequence}).
As explained in \cite{Bayer; Eisenbud 1995}, Theorem 1.6 the splittings of this induced extension of $\O_A$-modules
correspond bijectively to the desired extension $\varphi:B\ra C$ of the normalization map $A\ra C$ along the
inclusion $A\subset B$.
In other words,  we are looking for commutative diagrams
$$
\begin{CD}
\Omega^1_{B/k}\otimes\O_A@>>>\Omega^1_{A/k}\\
@Ad\varphi \otimes 1AA @AA A\\
\Omega^1_{C/k}\otimes\O_A@>>\id> \phantom{.}\Omega^1_{C/k}\otimes\O_A
\end{CD}
$$
of $\O_A$-modules. We thus happily arrived at a   linearization of the problem.

To proceed, we merely have to compute the sheaf of differentials
$\Omega^1_{B/k}$ and $\Omega^1_{C/k}$, together with  their restrictions to $A$.
We start with the nonreduced quadric $B$.
The $\O_B$-module $\Omega^1_{B/k}$ is freely generated by the differentials
$$
du,d\epsilon\quadand d(u^{-1}),d(\epsilon u^{-1})
$$
over the two open subsets $B\cap U$ and $B\cap U'$, respectively.
The corresponding 1-cocycle for the  locally free $\O_B$-module $\Omega^1_{B/k}$ is 
the $2\times 2$-matrix
$$
\begin{pmatrix}
u^{-2} & \epsilon u^{-2}\\
0 & u^{-1}
\end{pmatrix},
$$
because $d(u^{-1})=u^{-2}du$ on the overlap, and similarly for $d(\epsilon u^{-1})$.

We next turn to the rational cuspidal curve $C$.
The $\O_C$-module $\Omega^1_{C/k}$ is generated by the differentials
$$
d(u^2), d(u^3)\quadand d(u^{-1})
$$
over the two open subsets, respectively. On the first open subset,
we have a single relation $u^4d(u^2)=0$, because we are in characteristic $p=2$.
It follows that $\Omega^1_{C/k}$ modulo torsion
is invertible, with generators $d(u^3)$ and $d(u^{-1})$, and corresponding
1-cocycle $u^{-4}$. We refer to \cite{Schroeer 2005}, Section 3 for further results.

Now to the desired $\O_A$-linear map $d\varphi\otimes 1:\Omega^1_{C/k}\otimes\O_A\ra\Omega^1_{B/k}\otimes\O_A$.
Any such map is of the form
$$
d(u^3)\longmapsto u^2du+P(u)d\epsilon \quadand d(u^{-1})\longmapsto d(u^{-1})+Q(u^{-1})d(\epsilon u^{-1})
$$
for some polynomials $P(u)$ and  $Q(u^{-1})$ with coefficients from $k$. The 1-cocycles computed in the preceding paragraph
impose the condition $P(u)u^{-3}=Q(u^{-1})$. The upshot is that
the polynomials $P=P(u)$ of degree $\leq 3$ correspond to such $\O_A$-linear maps.
The corresponding morphism $\varphi:B\ra C$ is given in coordinates by
$$
u^2\longmapsto u^2,\quad u^3\longmapsto u^3+\epsilon P \quadand u^{-1}\longmapsto u^{-1}+\epsilon u^{-4}P.
$$
Throughout, we call $P$ the \emph{glueing polynomial} and write  it as
$$
P=\alpha_3u^3+\alpha_2u^2+\alpha_1u +\alpha_0,
$$
with scalars $\alpha_3,\ldots,\alpha_0\in k$.
For the constructions we have in mind it is important that the morphism $\varphi:B\ra C$ 
of degree two is flat. This condition depends on the constant term of the glueing polynomial:

\begin{proposition}
\mylabel{flat nonzero}
The morphism $\varphi:B\ra C$ is flat if and only if the constant term $\alpha_0$
in the glueing polynomial $P\in k[u]$ is nonzero.
\end{proposition}

\proof
Clearly, $\varphi$ is flat outside the singular point $c\in C$. 
So $\varphi$ is flat if and only if the Artin scheme $\varphi^{-1}(c)\subset B$ has length two.
Clearly, the fiber in question is the spectrum of the Artin ring
$$
k[u,\epsilon]/(u^2,u^3+\epsilon P) = k\oplus ku\oplus k\epsilon\oplus ku\epsilon/(\epsilon P).
$$
If $\alpha_0=0$, this is a $k$-vector space of dimension $d=3$ or $d=4$.
If $\alpha_0\neq 0$,   the residue class of $P$ is a unit, and the $k$-vector space
has dimension $d=2$.
\qed

\medskip
From now on we assume that the glueing polynomial 
$P=\alpha_3u^3+\ldots+\alpha_0$ has a nonzero constant term,
such that our morphism $\varphi:B\ra C$ is flat. 
Moreover, we regard our rational cuspidal curve as
$$
C=\Spec k[u^2,u^3+\epsilon P]\;\cup\; \Spec k[u^{-1}+\epsilon u^{-4}P].
$$
In other words, we view $\O_C$ as a subsheaf of $\O_B$ with respect to our morphism
$\varphi:B\ra C$, and not merely as a subsheaf of $\O_A$. 
We call the rational point $y_\infty\in C$ that constitutes the complement of the affine open subset
$\Spec k[u^2,u^3+\epsilon P]\subset C$ the \emph{point at infinity}. We shall
see that it plays a special role. This is already apparent in the following   fact.

\begin{proposition}
\mylabel{tschirnhausen invertible}
The $\O_C$-module $\shT=\varphi_*\O_B/\O_C$ is invertible, and isomorphic
to $\O_C(y_\infty)$, where  $y_\infty\in C$ is the point at infinity.
\end{proposition}

\proof
Since $\varphi:B\ra C$ is flat of degree two, $\shT$ must be invertible. We compute
$$
\deg(\shT)=\chi(\shT)-\chi(\O_C)=\chi(\O_B)-2\chi(\O_C)=1.
$$
Hence $H^0(C,\shT)$ is 1-dimensional. To compute a nonzero section, we
use the affine fpqc-covering $V\amalg V'\ra C$ given by the formal completion $V=\Spec k[[u^2,u^3+\epsilon P]]$
and the affine open subset $V'=\Spec k[u^{-1} + \epsilon u^{-4}P]$.
One easily sees that the residue classes of the unit $u\in\Gamma(V,\O_B)$ and the nilpotent $\epsilon u^{-1}P\in\Gamma(V',\O_B)$
generate the quotient sheaf $\shT$ as an $\O_C$-module. On the overlap $V\times_C V'$, we have
$(u^3+\epsilon P) u^{-2} \cdot u \cong \epsilon u^{-1}P$ modulo $\O_C$.
Using  
$$
\frac{1}{(u^3+\epsilon P) u^{-2}}=u^{-1}+\epsilon u^{-4} P,
$$
we see
that the  local sections $1\cdot u\in\Gamma(V,\shT)$ and
$(u^{-1}+\epsilon u^{-4} P)\cdot\epsilon u^{-1}P\in\Gamma(V',\shT)$ glue together and define a global section,
which vanishes  precisely at the point at infinity $y_\infty\in C$.
\qed

\section{The geometric construction}
\mylabel{The geometric construction}

We keep the notation from the preceding section, and 
use the flat morphism $\varphi:B\ra C$ to form the cocartesian square
\begin{equation}
\label{pushout diagram}
\begin{CD}
B @>>> X\\
@V\varphi VV @VV\nu V\\
C @>>> \phantom{.}Y.\\
\end{CD}
\end{equation}
The surface $Y$ is our desired del Pezzo surface, as we shall see
in due course.
Pushouts like the above exists as  
algebraic spaces according to a very general criterion of Artin \cite{Artin 1970}, Theorem 6.1.
The morphism $\nu:X\ra Y$ is the normalization, and the cartesian square is
also cocartesian. In particular, we have $\nu^{-1}(C)=B$.
Fortunately, we may immediately forget about the category of
algebraic spaces.

\begin{proposition}
\mylabel{space scheme}
The algebraic space $Y$ is a projective scheme.
\end{proposition}

\proof
It suffices to find an ample invertible
$\O_Y$-module.
As explained in \cite{Schroeer; Siebert 2002}, Proposition 4.1,
the pushout diagram yields an exact sequence of abelian sheaves
$$
1\lra\O_Y^\times\lra\O_X^\times\times\O_C^\times\lra\O_B^\times\lra 1,
$$
which results in an exact sequence of abelian groups
\begin{equation}
\label{pic sequence}
0\lra\Pic(Y)\lra\Pic(X)\oplus\Pic(C)\lra\Pic(B).
\end{equation}
Clearly, the preimage map $\Pic(C)\ra\Pic(B)$
is surjective. Hence there is an invertible $\O_C$-module $\shL$ with $\shL_B=\O_B(1)$.
Consequently, there is an invertible $\O_Y$-module $\O_Y(1)$ whose preimage on
$X$ is isomorphic to $\O_X(1)$. According to \cite{EGA IIIa}, Proposition 2.6.2,
the invertible sheaf $\O_Y(1)$ must be ample,
so the algebraic space $Y$ is projective.
\qed

\begin{remark}
The proof works under fairly general assumptions: $X$ might be any projective scheme,
$B$ a one-dimensional subscheme, and $\varphi:B\ra C$ a  morphism of curves that is generically
an isomorphism. 
\end{remark}

\medskip
It is not difficult to write down   the coordinate rings for the affine open covering $Y=V\cup V'\cup V''$
corresponding to the affine open covering $X=U\cup U'\cup U''$ defined in (\ref{covering X}).
Indeed, the diagram
\begin{equation}
\label{pullback diagram}
\begin{CD}
\O_B @<<< \O_X\\
@A\varphi AA @AA\nu A\\
\O_C @<<< \phantom{.}\O_Y\\
\end{CD}
\end{equation}
is cartesian, and this implies that
\begin{equation}
\label{covering Y}
\begin{gathered}
V=\Spec k[u^2,u^3+vP,v^2,v^2u,v^3,v^3u],\\
V'=\Spec k[u^{-2},u^{-1}+vu^{-4}P,v^2u^{-2},v^2u^{-3},v^3u^{-3},v^3u^{-4}],\\
V''=\Spec k[uv^{-1},v^{-1}].
\end{gathered}
\end{equation}
Let me explain this for the first open subset $V$:
Clearly, the six given elements $u^2,u^3+vP,\ldots, v^3u$ lie
in $\O_V$. Moreover, any monomial of the form $u^mv^n$ with $m,n\geq 2$
 is a monomial in the elements $u^2,v^2,v^2u,v^3,v^3u$.
Finally, the residue classes of $u^2,u^3+vP$ generate the quotient sheaf $\O_C|_V$.
Whence the given elements generate $\O_V$.

The  two open subsets $V,V'\subset Y$ need uncomfortably many generators.
It is possible to compute, with computer algebra,  a Gr\"obner basis for the ideal of relations,
but this sheds little light on the situation.
However, we shall see that things clear up
under passing to suitable localizations or completions.

The scheme $Y$ has as reduced singular locus the rational cuspidal curve $C\subset Y$. 
Our ultimate goal is to construct
twisted forms of $Y$ that are regular, or at least normal.
This can only happen if the singularities
on $Y$ are not too bad.
Recall that a locally noetherian scheme $S$ is called \emph{locally of complete intersection} if
for all points $s\in S$, the formal completion of the stalk $\O_{S,s}$ is of the
form $\O_{S,s}^\wedge=R/I$, where $R$ is a regular local noetherian ring  and $I\subset R$
is an ideal generated by a regular sequence (\cite{EGA IVd}, Definition 19.3.1).
The whole paper hinges on the following observation.

\begin{theorem}
\mylabel{complete intersection}
The  scheme $Y$ is locally of complete intersection.
\end{theorem}

\proof
We shall determine the local generators and relations explicitly.
This will be   useful later, when we  compute the   cotangent sheaf of the singular scheme $Y$.
To start with, consider the affine open subset $V\subset Y$ occurring in (\ref{covering Y}), with
coordinate ring $A=k[u^2,u^3+vP,v^2,v^2u,v^3,v^3u]$.
To simplify notation, we give names to the generators:
\begin{equation}
\label{abcec'f}
a=u^2,\quad b=u^3+vP,\quad c=v^2,\quad e=v^2u,\quad c'=v^3,\quad f=v^3u.
\end{equation}
The idea now is to localize so that fewer than six generators suffice.
First, let us take the affine open subset $V_{P^2}\subset V$ obtained by inverting the element
$$
P^2=\alpha_3^2a^3+\alpha_2^2a^2+\alpha_1^2a+\alpha_0^2\in A.
$$
Recall that $P=\alpha_3u^3+\ldots+\alpha_0$ is the glueing  polynomial in the indeterminate $u$
defining our glueing map $\varphi:B\ra C$. 
I contend that $A_{P^2}=k[a,b,e]_{P^2}$ as subrings inside the function field $k(u,v)$. Indeed, we have
$c=(b^2+a^3)/P^2$ and compute
$$
\begin{gathered}
cb=ea + c'(\alpha_0+\alpha_2a) +f(\alpha_1+\alpha_3a),\\
eb=ca^2+c'(\alpha_1a+\alpha_3a^2)+f(\alpha_0+\alpha_2a).
\end{gathered}
$$
The matrix of  coefficients at $c',f$ has determinant
$$
\det
\begin{pmatrix}
\alpha_0+\alpha_2a & \alpha_1+\alpha_3a\\
\alpha_1a+\alpha_3a^2 & \alpha_0+\alpha_2a
\end{pmatrix}
=P^2,
$$
and hence we may express $c',f$ in terms of $a,b,e, 1/P^2$. The upshot is that the canonical inclusion
$k[a,b,e]_{P^2}\subset A_{P^2}$ is bijective.

According to the work of Avramov \cite{Avramov 1975}, the property of being locally of complete
intersection is stable under localization. Therefore, it remains to check that 
$k[a,b,e]$ is of complete intersection. But this is   trivial:
We write the 2-dimensional ring $k[a,b,e]$ as a quotient of a polynomial ring in three indeterminates, and since the
latter is factorial, the ideal of relations is generated by a single element.
For later use, I write down the such a relation; it is
\begin{equation}
\label{first relation}
P^4e^2+b^4a+a^7.
\end{equation}

To continue, let us look at the affine open subset $V_Q\subset V$ given by localizing the
element $Q=a^2+c(\alpha_1+\alpha_3a)^2\in A$.
With this choice, I claim that $A_Q=k[a,b,c,c']_Q$. The argument is very similar to the one in the preceding
paragraph, and reveals how to come up with a denominator  like $Q$: We compute 
$$
\begin{gathered}
cb=ea+f(\alpha_1+\alpha_3a)+c'(\alpha_0+\alpha_2a),\\
c'b=ec(\alpha_1+\alpha_3a)+fa+c^2(\alpha_0+\alpha_2a).
\end{gathered}
$$
The coefficients at $e,f$ comprise a matrix, whose determinant is
$$
\det
\begin{pmatrix}
a & \alpha_1+\alpha_3a\\
c(\alpha_1+\alpha_3a) & a
\end{pmatrix}
=
a^2+c(\alpha_1+\alpha_3a)^2 = Q.
$$
Hence we may express the generators $e,f$ in terms of $a,b,c,c',1/Q$, and therefore $A_Q=k[a,b,c,c']_Q$.
It remains to see that $R=k[a,b,c,c']$ is a complete intersection.
The generators $a,c\in R$ are algebraically independent, and we have relations
\begin{equation}
\label{second relation}
b^2+a^3+cP^2\quadand c'^2=c^3.
\end{equation}
I claim that the canonical surjection
$$
R'=k[a,c][x_1,x_2]/(x_1^2+a^3+cP^2, x_2^2+c^3)\lra R,\quad x_1\longmapsto b,\quad x_2\longmapsto c'
$$
is bijective, where $x_1,x_2$ are indeterminates. Indeed, both rings in question are 2-dimen\-si\-onal, and $R$ is integral,
so it suffices to check that the ring $R'$ on the left 
is integral. The inclusion $k[a,c]\subset R'$ is purely inseparable, so $\Spec(R')$ is irreducible.
Being a complete intersection, the affine scheme $\Spec(R')$ has no embedded components.
Hence, to check that the 2-dimensional ring $R'$ is reduced we may replace it by the
local Artin ring $k(a,c)[x_1,x_2]/(x_1^2+a^3+cP^2, x_2^2+c^3)$.
By Lemma \ref{Artin field} below, it suffices to check that
the differentials 
$$
d(a^3+cP^2)=(a^2+\alpha_1^2c +\alpha_3^2ca^2)da +P^2dc \quadand d(c^3)=c^2dc
$$
from $\Omega^1_{k(a,c)/k}$ are linearly independent, and the latter is obvious.
If follows that the affine open subset $V\subset Y$ is locally of complete intersection,
because the open subsets $V_{P^2}, V_Q\subset V$ cover the   singular locus $\Sing(V)=V\cap C$.

To finish the proof, it suffices to see that the formal completion  
$$
R=k[[u^{-2},u^{-1}+vu^{-4}P,v^2u^{-2},v^2u^{-3},v^3u^{-3},v^3u^{-4}]]
$$
of the second affine open subset $V'\subset Y$ is  a complete intersection.
We introduce names
\begin{equation}
\label{xyz}
x=u^{-1}+vu^{-4}P \quadand y=v^2u^{-2} \quadand z=v^3u^{-3}.
\end{equation}
The equation
$
x^2= u^{-2} + y(\alpha_3^2+\alpha_2^2u^{-2}+\alpha_1^2xu^{-4} +\alpha_0^2u^{-6})
$,
viewed as a recursion relation for $u^{-2}$, reveals that the generator $u^{-2}$
is a formal power series in $x^2,y$.
Next, we decompose the glueing polynomial $P=\alpha_3u^3+\ldots+\alpha_0$ into
even and odd part $P=P_\text{ev}+P_\text{odd}$, such that
$P_\text{ev}u^{-2}$ and $P_\text{odd}u^{-3}$ are  polynomials in $u^{-2}$.
Computing
$$
\begin{gathered}
xy= v^2u^{-3} + v^3u^{-4} P_\text{ev}u^{-2} + zP_\text{odd}u^{-3},\\
xz=v^2u^{-3} yP_\text{ev}u^{-2} + v^3u^{-4} + y^2P_\text{odd}u^{-3},
\end{gathered}
$$
we see that the matrix of coefficients at $v^2u^{-3},v^3,u^{-4}$ has determinant
$$
\det\begin{pmatrix}
1 & P_\text{ev}u^{-2}\\
yP_\text{ev}u^{-2} & 1
\end{pmatrix}
=1 - yP_\text{ev}^2u^{-4},
$$
which is a unit. Whence it is possible to express the generators $v^2u^{-3},v^3u^{-4}$
as formal power series
in $x,y,z$, so the inclusion $k[[x,y,z]]\subset R$ is bijective.
But the ring $k[[x,y,z]]$ is obviously a complete intersection,
with relation $y^3+z^2=0$.
\qed

\medskip
In the preceding proof, we needed the following fact.

\begin{lemma}
\mylabel{Artin field}
Let $K$ be a field of characteristic $p>0$,
and $f_1,\ldots,f_n\in K$ elements so that
the differentials $df_1,\ldots, df_n\in\Omega^1_{K/K^p}$
are linearly independent.
Then the local Artin ring 
$K[x_1,\ldots,x_n]/(x_1^p-f_1,\ldots,x^p_n-f_n)$
is a field.
\end{lemma}

\proof
Let $A=K[x_1,\ldots,x_n]$ be the polynomial algebra
and set $y_i=x_i^p-f_i$.
Let $\maxid\subset A$ be the maximal ideal containing the
ideal $(y_1,\ldots,y_n)$. 
It suffices to check that the residue classes $y_i\in \maxid/\maxid^2$
are linearly independent, according to \cite{EGA IVa}, Proposition 17.1.7.

By assumption, we find $K^p$-derivations $D_i:K\ra K$
with $D_i(f_j)=\delta_{ij}$ (Kronecker delta).
The cotangent sequence
$$
0\lra\Omega^1_{K/K^p}\otimes A\lra\Omega^1_{A/K^p}\lra\Omega^1_{A/K}\lra 0
$$
is exact and splits, because the ring extension $K\subset A$ is smooth.
Hence we may extend our $D_i$ to $K^p$-derivations $D_i:A\ra A$,
which have $D_i(y_j)=\delta_{ij}$.
These derivations induce linear maps $D_i:\maxid/\maxid^2\ra A/\maxid$.
If follows that the residue classes of $y_i$ are linearly independent.
\qed

\begin{remark}
The condition that differentials $df_1,\ldots,df_n\in\Omega^1_{K/K^p}$ are linearly independent
exactly means that the elements $f_1,\ldots, f_n\in K$ are \emph{$p$-linearly independent}.
\end{remark}

\section{Picard scheme and dualizing sheaf}
\mylabel{Picard scheme and dualizing sheaf}

We keep the notation from the preceding section, such that
$Y$ is an integral projective surface  locally of complete intersection,
which is defined by the cocartesian square (\ref{pushout diagram}).
It is nonnormal, with reduced nonsmooth locus $C\subset Y$ and 
normalization $X=\PP^2$. In this section  we study invertible sheaves
on $Y$. In some sense, everything reduces to the curve $C\subset Y$:

\begin{proposition}
\mylabel{pic Y}
The restriction map $\Pic_{Y/k}\ra\Pic_{C/k}$ of Picard schemes
is an isomorphism.
\end{proposition}

\proof
The exact sequence of abelian sheaves
$1\ra\O_Y^\times\ra\O_X^\times\times\O_C^\times\ra\O_B^\times\ra 1$
induces an exact sequence of Picard schemes
$$
0\lra\Pic_{Y/k}\lra\Pic_{X/k}\times\Pic_{C/k}\lra\Pic_{B/k}.
$$
Hence it suffices to check that the restriction map $\Pic_{X/k}\ra\Pic_{B/k}$
is an isomorphism. Recall that $B$ is the first order infinitesimal neighborhood
of a line $A$ inside $X=\PP^2$.
Clearly, the restriction map $\Pic_{X/k}\ra\Pic_{A/k}$ is an isomorphism, so it suffices
to check that the restriction map $\Pic_{B/k}\ra\Pic_{A/k}$ is an isomorphism.
Let $\shI\subset\O_B$ be the ideal for the closed embedding $A\subset B$.
Then $\shI^2=0$ and  
we have an exact sequence of abelian sheaves $0\ra\shI\ra\O_B^\times\ra\O_A^\times\ra 1$,
which gives an exact sequence
$$
H^1(B,\shI)\lra \Pic_{B/k}\lra\Pic_{A/k}\lra H^2(B,\shI).
$$
The outer terms vanish, because the abelian sheaf $\shI$ is isomorphic to the  $\O_A$-module $\O_A(-1)$,
whence the assertion holds.
\qed

\medskip
We conclude that the Picard scheme $\Pic_{Y/k}=\Pic_{C/k}$ is reduced and therefore smooth, and 1-dimensional.
In particular, its tangent space $H^1(Y,\O_Y)$ is 1-dimensional.
The Picard scheme sits inside a split extension
\begin{equation}
\label{picard extension}
0\lra\GG_a\lra\Pic_{Y/k}\lra\ZZ\lra 0.
\end{equation}
The map on the right is given by sending an invertible $\O_Y$-module $\shL$
to the degree of the restriction  $\shL_C$. To simplify notation, we set $\deg(\shL)=\deg(\shL_C)$ and
call this integer the \emph{degree} of $\shL$. 

\begin{proposition}
\mylabel{ample degree}
An invertible $\O_Y$-module $\shL$ is ample if and only if
$\deg(\shL)>0$.
\end{proposition}

\proof
According to  \cite{EGA IIIa}, Proposition 2.6.2, the invertible sheaf $\shL$ is
ample if and only if its preimage $\nu^*(\shL)$ is ample.
If follows easily from the definition of $d=\deg(\shL)$
that $\nu^*(\shL)=\O_X(d)$.
Hence $\shL$ is ample if and only if $d>0$.
\qed

\medskip
Being locally of complete intersection, the proper scheme
$Y$ also has an invertible dualizing sheaf $\omega_Y$.
It is straightforward to compute its degree:

\begin{proposition}
\mylabel{preimage dualizing}
The degree of the dualizing sheaf is $\deg(\omega_Y)=-1$.
\end{proposition}

\proof
Consider the Tschirnhausen module $\shT=\nu_*\O_X/\O_Y$.
Its annulator ideal $\cid\subset\O_Y$ is called the 
\emph{conductor ideal} for the normalization map $\nu:X\ra Y$.
According to Proposition \ref{tschirnhausen invertible}, the $\O_Y$-module $\shT$ is
an invertible $\O_C$-module, and hence $\cid=\O_Y(-C)$.
Since the   square (\ref{pushout diagram}) defining $Y$ is cocartesian,
the induced ideal on $X$ satisfies $\cid=\cid\O_X$.
As the square is also cartesian, we have $\cid\O_X=\O_X(-B)$.

The conductor   is closely related to duality:
The equality $\cid=\shHom(\nu_*\O_X,\O_Y)$ shows that
the conductor ideal has a natural $\O_X$-module structure, and
coincides with the relative dualizing sheaf $\omega_{X/Y}$.
The latter satisfies $\omega_X=\omega_{X/Y}\otimes\nu^*\omega_Y$.
Clearly, the projective plane $X=\PP^2$ has dualizing sheaf
$\omega_X=\O_X(-3)$. Together with $\omega_{X/Y}=\O_X(-2)$,
it follows $\nu^*(\omega_Y)=\O_X(-1)$.
\qed

\medskip
In particular, the dualizing sheaf $\omega_Y$ is antiample.
We thus call our $Y$ a \emph{nonnormal del Pezzo surface}.
Let me point out that its irregularity is $h^1(\O_Y)=1$,
which is highly unusual for del Pezzo surfaces, even for singular ones
(compare \cite{Hidaka; Watanabe 1981}, Corollary 2.5 and \cite{Schroeer 2001}, Theorem 2.2 and \cite{Reid 1994}).

To determine the isomorphism  class of $\omega_Y$ in the Picard group, it
suffices to compute its restriction to $C$.
Recall that the Tschirnhausen module $\shT=\nu_*\O_X/\O_Y$
is the invertible $\O_C$-module $\O_C(y_\infty)$, where
$y_\infty\in C$ is the point at infinity.

\begin{proposition}
\mylabel{restriction dualizing}
With the preceding notation, we have $\omega_Y|_C=\O_C(-y_\infty)$.
\end{proposition}

\proof
We consider the relative dualizing sheaf $\omega_{C/Y}=\shExt^1_{\O_Y}(\O_C,\O_Y)$,
which satisfies $\omega_C=\omega_{C/Y}\otimes\omega_Y|_C$.
The exact sequence
$0\ra\cid\ra\O_C\ra\O_C\ra0$ yields an exact sequence
$$
0\lra\shHom_{\O_Y}(\O_Y,\O_Y)\lra\shHom_{\O_Y}(\cid,\O_Y)\lra\shExt^1_{\O_Y}(\O_C,\O_Y)\lra 0.
$$
We have an obvious inclusion $\O_X\subset\shHom(\cid,\O_Y)$, and we now
check that it is bijective.
The composition map
$$
\shHom_{\O_Y}(\O_X,\O_Y)\otimes_{\O_X} \shHom_{\O_X}(\cid,\O_X)\lra\shHom_{\O_Y}(\cid,\O_Y)
$$
is bijective, because the conductor ideal is invertible as $\O_X$-module.
For the same reason, the  evaluation map 
$$
\shHom_{\O_X}(\cid,\O_X)\otimes_{\O_X}\cid\ra\O_X
$$
is bijective. 
Composing the previous maps, we obtain 
a chain of inclusion $\O_X\subset \shHom(\cid,\O_Y)\subset\O_X$, which clearly is bijective.
The upshot is that the relative dualizing sheaf $\omega_{C/Y}$ coincides
with the Tschirnhausen module $\shT$. Using that $\omega_C=\O_C$ 
and Proposition \ref{tschirnhausen invertible}, we deduce the assertion.
\qed

\section{Cartier divisors and Weil divisors}
\mylabel{Cartier divisors and Weil divisors}

Our del Pezzo surface $Y$ has a natural polarization furnished by the
ample invertible sheaf $\omega_Y^\vee$.
Given any Weil divisor $D$ on $Y$, we define its \emph{degree}
by the intersection number $\deg(D)=\omega_Y^\vee\cdot D$.
In this section we have a closer look at curves $D\subset Y$ of degree one.
Much of the geometry of $Y$ is captured by these curves.
Any such curve is the schematic image of a line $L$ on the normalization
$X=\PP^2$. To begin with, we compute some  
cohomology groups.

\begin{proposition}
\mylabel{euler characteristic}
Let $\shL$ be an invertible $\O_Y$-module of degree $d$.
Then we have $\chi(\shL)=d(d+1)/2$. Moreover, $H^2(Y,\shL)=0$ for $d\geq 0$,
and $H^1(Y,\shL)=0 $ for $d\geq 1$.
\end{proposition}

\proof
The exact sequence
$0\ra\shL\ra\shL_X\ra\shL\otimes\shT\ra 0$
of coherent sheaves gives $\chi(\shL)=\chi(\shL_X)-\chi(\shL\otimes\shT)$.
Here $\shT=\O_X/\O_Y$ is the Tschirnhausen module.
We clearly have $\chi(\shL_X)=(d+2)(d+1)/2$.
According to Proposition \ref{tschirnhausen invertible},   $\shT$ is an invertible $\O_C$-module of degree one.
It follows that $\chi(\shL\otimes\shT)=d+1$, and the assertion on the
Euler characteristic follows.

Now suppose $d\geq 0$. Then the term on the left in the exact sequence
$$
H^1(C,\shL\otimes\shT)\lra H^2(Y,\shL)\lra H^2(X,\shL_X)
$$
vanishes. The term on the right is Serre dual to $H^0(X,\shL_X^\vee(-3))$,
and vanishes as well. Hence $H^2(Y,\shL)=0$.

Finally, suppose that the degree is $d\geq 1$.
We now use the short exact sequence 
$0\ra\shL\ra\shL_X\oplus\shL_C\ra\shL_B\ra 0$, which gives an exact sequence
$$
H^0(X,\shL_X)\oplus H^0(C,\shL_C)\lra H^0(B,\shL_B)\lra H^1(Y,\shL)\lra H^1(X,\shL_X)\oplus H^1(C,\shL_C).
$$
The sum on the right vanishes, and the map on the left is surjective.
The latter follows from the exact sequence
$$
H^0(X,\shL)\lra H^0(B,\shL_B)\lra H^1(X,\shL(-2))=0.
$$
The upshot is that $H^1(Y,\shL)=0$.
\qed

\medskip
For $d=1$, this means that $H^0(Y,\shL)$ is 1-dimensional.
In other words, each invertible sheaf $\shL$ of degree one, 
there is precisely one effective Cartier divisor $D\subset Y$
with $\shL=\O_Y(D)$. Note that this applies in particular to the antidualizing
sheaf $\shL=\omega_Y^\vee$.  It turns out that the position 
of these Cartier divisors is determined by restricting to
the reduced singular locus $C\subset Y$:

\begin{proposition}
\mylabel{restriction bijective}
Let $\shL$ be an invertible sheaf of degree $d=1$.
Then the restriction map  $H^0(Y,\shL)\ra H^0(C,\shL_C)$ 
is bijective.
 \end{proposition}

\proof
The exact sequence $0\ra\shL_X(-B)\ra\shL_X\ra\shL_B\ra 0$ gives
an exact sequence
$$
H^0(X,\shL_X(-B))\lra H^0(X,\shL_X)\lra H^0(B,\shL_B)\lra H^1(X,\shL_X(-B)).
$$
Both outer terms vanish, and hence the restriction map 
$H^0(X,\shL_X)\ra H^0(B,\shL_B)$ is bijective.

Using  the exact sequence $0\ra\shL\ra\shL_X\oplus\shL_C\ra\shL_B\ra 0$,
we obtain a short exact sequence
$$
0\lra H^0(Y,\shL)\lra H^0(X,\shL_X)\oplus H^0(C,\shL_C)\lra H^0(B,\shL_B)\lra 0.
$$
In light of the preceding paragraph, the restriction map
$H^0(Y,\shL)\ra H^0(C,\shL_C)$ must be bijective.
\qed

\medskip
We conclude that given a rational point $y\in C$ in the smooth locus of $C$,
there is precisely one Cartier divisor $D\subset Y$ of degree one passing through
$y$.
Of course, there is a continuous family of Weil divisors of degree one passing through that point.
Any such Weil divisor is the image of a unique line on $X=\PP^2$.
How to distinguish between  such Cartier divisors and Weil divisors?

\begin{proposition}
\mylabel{cartier divisor}
Let $L\subset X$ be a line not contained in the conductor locus $B\subset X$,
and  $D\subset Y$ be its image, and $y\in C\cap D$ be the unique
intersection point with the singular locus $C\subset Y$.
Then the Weil divisor $D$ is Cartier if and only if
the schematic preimage $\nu^{-1}(y)\subset X$ is contained in $L$.
\end{proposition}

\proof
Suppose that $D\subset Y$ is Cartier.
Then the preimage $\nu^{-1}(D)\subset X$ is Cartier as well, and clearly contains $\nu^{-1}(y)$.
The inclusion $L\subset\nu^{-1}(D)$ is an equality outside the conductor locus $B\subset X$.
Since both subschemes are Cartier, they must be equal. If follows that $\nu^{-1}(y)\subset L$.

Now suppose $\nu^{-1}(y)\subset L$. Let $D'\subset Y$ be the unique
Cartier divisor of degree one passing through $y$. Its preimage
$L'=\nu^{-1}(D')\subset X$ is a line containing $\nu^{-1}(y)$, which
is an Artin subscheme of length two. However, through any Artin subscheme of length two on $X=\PP^2$, there
passes precisely one line. We conclude $L=L'$, whence $D=D'$ is Cartier.
\qed

\medskip
We finally determine what kind of scheme a Weil divisor
of degree one is. 

\begin{proposition}
\mylabel{cartier weil}
Let $D\subset Y$ be a Weil divisor of degree one.
\renewcommand{\labelenumi}{(\roman{enumi})}
\begin{enumerate}
\item
If $D\subset Y$ is Cartier or if $D=C$, then the curve $D$ is isomorphic to the rational cuspidal curve
with arithmetic genus $p_a=1$.
\item
If $D\subset Y$ is not Cartier and $D\neq C$, then the curve $D$ is isomorphic to the projective line $\PP^1$.
\end{enumerate}
\end{proposition}

\proof
The case $D=C$ is clear, because $C$ is by definition  the rational cuspidal curve with $p_a=1$.
Now suppose that $D$ is Cartier.
Then $\omega_D=\omega_Y(D)\mid_D$ has degree zero. It follows that $-2\chi(\O_C)=\deg(\omega_D)=0$,
whence $h^1(\O_C)=1$.  Clearly, $D$ is birational and homeomorphic to the projective line,
and the assertion follows.

Finally, suppose that $D\neq C$ is not Cartier. Let $L\subset X$ be the unique line
with $D=\nu(L)$, and consider the birational morphism $f:L\ra D$.
According to Proposition \ref{cartier divisor}, the fiber $f^{-1}(y)=L\cap\nu^{-1}(y)$
is an Artin scheme of length one. It follows that $f$ is an isomorphism.
\qed

\medskip
To close this section, we look again at the antidualizing sheaf $\omega_Y^\vee$.
We know from  Proposition \ref{euler characteristic} that there is \emph{only one} effective anticanonical divisor $D\subset Y$.
We note in passing  an interesting consequence: The anticanonical divisor $D\subset Y$ is invariant under any automorphism
of $Y$.

\section{The tangent sheaf}
\mylabel{The tangent sheaf}

We keep the notation from the preceding sections, such that $Y$ is a nonnormal
del Pezzo surface.
To construct twisted forms of $Y$, we have to understand
the group scheme $\Aut_{Y/k}$ and its Lie algebra $H^0(Y,\Theta_{Y/k})$.
In this section we shall see that the tangent sheaf $\Theta_{Y/k}$ is locally free
of rank two, and that it is not difficult to determine its global sections.
This is somewhat surprising, because the cotangent sheaf $\Omega^1_{Y/k}$
is not locally free along the singular curve $C\subset Y$.
However, we shall see that the trouble only comes from the
torsion subsheaf $\tau\subset \Omega^1_{Y/k}$. This effect seems
to be special to positive characteristics.

\begin{theorem}
\mylabel{cotangent sheaf}
The coherent $\O_Y$-module $\Omega^1_{Y/k}/\tau$ is locally free of rank two.
\end{theorem}

\proof
This is a local problem in $Y$. First, consider the affine open subset $V\subset Y$
that is the spectrum of
$$
k[u^2,u^3+vP, v^2,v^2u,v^3,v^3u]=k[a,b,c,e,c',f],
$$
as in Equation (\ref{abcec'f}).
As explained in the proof for Proposition \ref{complete intersection},
it is advisable to localize further, using $P^2=\alpha_3^2a^3+\alpha_2^2a^2+\alpha_1^2a+\alpha_0^2$ and $Q=a^2+c(\alpha_1+\alpha_3a)^2$
as denominators.
We saw that $V_{P^2}$ is an open subset inside the spectrum of
$$
k[a,b,e]/(P^4e^2+b^4a+a^7).
$$
The module of differentials is generated by $da,db,de$ modulo the
relation $(b^4+a^6)da$. Since $Y$ is generically smooth, the coefficient $b^4+a^6$ must be a regular element,
hence the differential $da$ is torsion.
We conclude that the differentials $db,de$ form a basis of $\Omega^1_{Y/k}$ modulo torsion
over $V_{P^2}\subset Y$.
A similar argument using Equation (\ref{second relation})  gives that the differentials $db,dc'$ form a basis
for $\Omega^1_{Y}$ modulo torsion over $V_Q\subset Y$.

It remains to treat the affine open subset $V'\subset Y$. Since
the point at infinity $y_\infty\in V'$ is the only singularity not contained in $V$,
it suffices to consider the formal completion
$R=k[[x,y,z]]/(y^3-z^2)$ of $\O_{Y,y_\infty}$ as in Equation (\ref{xyz}).
The separated completion $\widehat{\Omega}^1_{R/k}$ modulo torsion is free, with basis
$dx,dz$. Since the two functions $x=u^{-1}+vu^{-4}P$ and $z=u^3v^{-3}$ are already contained
in $\O_{Y,y_\infty}$, it follows that $dx,dz$ are a basis modulo torsion in
some affine neighborhood of $y_\infty\in Y$.
\qed

\begin{corollary}
\mylabel{tangent sheaf}
The tangent sheaf $\Theta_{Y/k}$ is   locally free of rank two.
\end{corollary}

\proof
Dualizing the exact sequence
$
0\ra\tau\ra\Omega^1_{Y/k}\ra\Omega^1_{Y/k}/\tau\ra 0
$,
we see that the canonical map $\shHom(\Omega^1_{Y/k}/\tau,\O_Y)\ra\shHom(\Omega^1_{Y/k},\O_Y)$ is bijective.
\qed

\medskip
Our next task is to compute the Lie algebra of global sections for the tangent sheaf.
We are mostly interested in the behavior of derivations near the singular locus $C\subset Y$,
whence we shall describe $ H^0(Y,\Theta_{Y/k})$ as a subalgebra of
$H^0(V_{P^2},\Theta_{Y/k})$.
We just saw that $\Omega^1_{Y/k}$ modulo torsion has basis a on $V_{P^2}\subset Y$ given by $db,de$.
We denote by $D_b,D_e$ the dual basis of $\Theta_{Y/k}$ on $V_{P^2}$.
The following result gives an implicit description of $H^0(Y,\Theta_{Y/k})$,
which will give enough information for our purposes.

\begin{proposition}
The Lie algebra $H^0(Y,\Theta_{Y/k})$
consists of  all derivations of the form $fD_b+gD_e$, where
$f,g\in k[u^2,u^3+vP, v^2,v^2u,v^3,v^3u]$ are polynomials so that the rational functions
\begin{equation}
\label{rational functions}
f\frac{u}{v^2P}+g\frac{Pv+P'uv+u^3}{v^4P}
\quadand
f\frac{1}{v^2P}+g\frac{u^2+P'v}{v^4P}
\end{equation}
are contained  in $k[u/v,1/v]$.
\end{proposition}

\proof
Since the cotangent sheaf satisfies Serre's condition $(S_2)$ and the
complement of $V_{P^2}\cup V''\subset Y$ is finite, the restriction map
$H^0(Y,\Theta_{Y/k})\ra H^0(V_{P^2}\cup V'',\Theta_{Y/k})$
must be bijective.
The latter group   is the kernel of the difference map
\begin{equation}
\label{cech cohomology}
H^0(V_{P^2},\Theta_{Y/k})\oplus H^0(V'',\Theta_{Y/k})\lra H^0(V_{P^2}\cap V'',\Theta_{Y/k})
\end{equation}
coming from \v{C}ech cohomology. To determine the kernel, we
first compute with differentials rather than derivations: 
$$
\begin{gathered}
db=(u^2 + P'v)du + Pdv,\quad  de=v^2du\\
d(u/v)=1/vdu+u/v^2dv,\quad d(1/v)=1/v^2dv,
\end{gathered}
$$
as follows from (\ref{abcec'f}). Consequently,
\begin{equation}
\label{base change}
A=\begin{pmatrix}
u^2 +P'v & v^2\\
P & 0
\end{pmatrix},\quad
B=\begin{pmatrix}
1/v & 0\\
u/v^2 & 1/v^2
\end{pmatrix}
\end{equation}
are the base change matrices for base changes from $db,de$ to $du,dv$, and from 
$d(u/v),d(1/v)$ to $du,dv$, respectively.
It follows that $B^{-1}A$ is the base change matrix from $db,de$ to $d(u/v),d(1/v)$,
and whence 
$$
{}^t\!(B^{-1}A)^{-1}=\begin{pmatrix}
u/v^2P & (Pv+P'uv+u^3)/v^4P\\
1/v^2P & (u^2 + P'v)/v^4P
\end{pmatrix}
$$
is the base change matrix from  the dual basis $D_b,D_e$ to $D_{u/v},D_{1/v}$.
Using this base change matrix, we compute the kernel in the exact sequence (\ref{cech cohomology}),
and the assertion follows.
\qed

\medskip
Recall that the normalization of $Y$ is the projective plane $X=\PP^2$.
Pulling back to $X$, we see that $H^0(X,\Theta_{Y/k}\otimes\O_X)$ is given
by derivations of the form $fD_b+gD_e$, where
$f,g\in k[u,v]$ are polynomials so that the two rational functions
in (\ref{rational functions})  lie in $k[u/v,1/v]$.
In particular, the rational derivation $\delta=PD_e$ defines a global section of $\Theta_{Y/k}\otimes\O_X$.
It does not, however, always come from a global section of $\Theta_{Y/k}$:

\begin{corollary}
Suppose the glueing polynomial $P$ is even, that is, $P=\alpha_2u^2+\alpha_0$.
Then the rational vector field $\delta=PD_e$ 
lies in $H^0(Y,\Theta_{Y/k})$. Moreover, we have $\delta\circ\delta=0$.
\end{corollary}

\proof
The first statement follows from the preceding proposition.
We have $\delta(e)=\alpha_2a+\alpha_0$ and $\delta(b)=0$, and therefore
$\delta\circ\delta=0$.
\qed

\medskip
From now on we assume that the glueing polynomial $P$ is even.
Using  the   base change matrices in (\ref{base change}), we easily express
the global vector field $\delta=PD_e$ in terms of other rational derivations, and obtain
$$
\delta = PD_e= (Pv+u^3)v^{-4}D_{u/v}+ u^2v^{-4}D_{1/v} = Pv^{-2}D_u +  u^2v^{-2}D_v.
$$
In the next section, we shall interpret $\delta$ as a group scheme action of $\alpha_2$ on $Y$.
The fixed points for the group scheme action correspond to the zeros of the vector field.
By definition, $\delta=PD_e$ has no zeros on the open subset $V_{P^2}$.
Clearly, it vanishes on $V''$ to first order along the closed subscheme given by $1/v=0$.
It remains to determine the precise behavior of the vector field near the
curve of singularities $C\subset Y$.

\begin{proposition}
\mylabel{zeros field}
The only zero of the global vector field $\delta=PD_e$ 
lying on the singular curve $C\subset Y$ is the point
at infinity $y_\infty\in C$.
With respect to the formal coordinates $\O_{Y,y_\infty}^\wedge=k[[x,y,z]]/(y^3-z^2)$,
we have $\delta= xD_z$.
\end{proposition}

\proof
Using $\delta =Pv^{-2}D_u +  u^2v^{-2}D_v$ and 
the definitions of $x,y,z$ in Equation (\ref{xyz}), one computes $\delta(x)=\delta(y)=0$ and $\delta(z)=x$,
whence $\delta=xD_z$.
It remains to write $\delta$ on $V_Q\subset V$ in terms of the basis $D_b,D_{c'}$, where $c'=v^3$.
Note that $Q=a+c(\alpha_3u^3+\alpha_1u)=a$, because we assume the glueing polynomial to be even.
We have $\delta(b)=0$ and compute $\delta(v^3)=u^2=a$, whence $\delta=aD_{c'}$ has no zero on
the affine open subset  $V_Q\subset Y$.
\qed

\section{Splitting type of tangent sheaf}
\mylabel{Splitting type of tangent sheaf}

In this section we   determine the restriction $\Theta_{Y/k}|_D$
of the tangent sheaf to Weil divisors $D\subset Y$ of degree one.
This will show that our choice of global vector field $\delta=PD_e$ is, in some sense,
the best possible choice.

The computation with the tangent sheaf is rather easy, because we may dually work with
the cotangent sheaf modulo torsion.
The latter sits in an exact sequence
\begin{equation}
\label{exact sequence}
\shI/\shI^2\lra\Omega^1_{Y/k}\otimes\O_D\lra\Omega^1_{D/k}\lra 0,
\end{equation}
where $\shI=\O_Y(-D)$ is the ideal of the Weil divisor.
For the following arguments, note that the torsion subsheaf $\tau\subset\Omega^1_{Y/k}$ is
locally a direct summand, because $\Omega^1_{Y/k}/\tau$ is locally free,
whence this subsheaf commutes with base change.

\begin{proposition}
\mylabel{four one}
Let $D\subset Y$ be a Cartier divisor of degree one.
Then the tangent sheaf splits as $\Theta_{Y/k}|_D\simeq \Theta_{D/k}\oplus\omega_Y^\vee|_D$.
Both summands are invertible $\O_D$-modules, of degree four and one,
respectively.
\end{proposition}

\proof
The $\O_D$-module $\shI/\shI^2$ is invertible of degree $-D^2=-1$, and the canonical
map $\shI/\shI^2\ra\Omega^1_{Y/k}\otimes\O_D$ on the left in
the cotangent sequence (\ref{exact sequence}) is injective, because
$D\subset Y$ is Cartier.
According to Proposition \ref{cartier weil}, the scheme $D$ is the rational cuspidal curve
with arithmetic genus $p_a=1$. It follows $\shI/\shI^2\simeq \omega_Y|_D$, and that $\Omega^1_{D/k}$ modulo torsion is
an invertible sheaf of degree $-4$.
Since $Y$ is smooth at the generic point $\eta\in D$, the torsion  in $\Omega^1_{Y/k}\otimes\O_D$
maps to the torsion of $ \Omega_{D/k}^1$, and this map must be bijective
because $\Omega^1_{Y/k}/\tau$ is locally free of rank two, which contains $\shI/\shI^2$
locally as a direct summand. 
The result now follows by taking duals, and the fact that there are no  nontrivial
extension of invertible sheaves  of degree one by invertible sheaves of degree four on $D$.
\qed

\medskip
In particular, we see that the invertible sheaf $\det(\Theta_{Y/k})$ has degree five.

\begin{proposition}
\mylabel{two three}
Let $D\subset Y$ be a Weil divisor of degree one that is not Cartier and not
the singular curve $C$. Then $\Theta_{Y/k}|_D\simeq \O_D(3)\oplus\O_D(2)$ is a direct
sum of invertible $\O_D$-modules of degree two and three.
\end{proposition}

\proof
In this case, Proposition \ref{cartier weil} tells us that $D\simeq\PP^1$.
Hence $\Omega^1_{D/k}$ is invertible of degree $-2$. 
It follows that the torsion of $\Omega^1_{Y/k}\otimes \O_D$ maps to zero
in $\Omega^1_{Y/k}$. The induced map $(\Omega_{Y/k}^1/\tau)\otimes\O_D\ra\Omega^1_{D/k}$
has invertible kernel, which must have  degree $-3$.
The result follows after dualizing, and the fact that there are non nontrivial extensions
of $\O_{\PP^1}(3)$ by $\O_{\PP^1}(2)$.
\qed

\begin{proposition}
\mylabel{four one again}
Let $C\subset Y$ be the reduced singular locus.
Then the tangent sheaf splits as $\Theta_{Y/k}|_C\simeq \Theta_{C/k}\oplus\O_C(y_\infty)$.
Both summands are invertible, of degree four and one, respectively.
\end{proposition}

\proof
I claim that the torsion in $\Omega^1_{Y/k}$ maps to the torsion in $\Omega^1_{C/k}$.
It suffices to check this on the formal completion $R'=k[[x,y,z]]/(y^3+z^2)$
of the affine open subset $V'\subset Y$ at the point at infinity.
The curve $C\subset Y$ has ideal $(y,z)$, and the torsion is generated by $dy$,
whence the claim follows.
Whence we have an exact sequence
$$
0\lra \shK\lra\Omega^1_{Y/k}\otimes\O_C/(\text{torsion})\lra \Omega^1_{C/k}/(\text{torsion})\lra 0
$$
for some coherent $\O_C$-module $\shK$.
Since both terms on the right are locally free, the $\O_C$-module $\shK$ must be invertible.
It must have degree $-1$, because $\Omega_{C/k}^1$ modulo torsion is invertible of degree $-4$.
Such an extension of $\O_C$-modules must split. Dualizing it, we obtain $\Theta_{Y/k}|_C=\Theta_{C/k}\oplus\shK^\vee$.

It remains to see that $\shK^\vee\simeq\O_C(y_\infty)$. We use our global vector field $\delta=PD_e$.
By Proposition \ref{zeros field}, the restriction $\delta\otimes 1$ vanishes only at $y_\infty\in C$,
and has vanishing order one there.
Decompose $\delta\otimes 1=\delta'+\delta''$, where $\delta'$ is a global vector field on $C$,
and $\delta''\in H^0(C,\shK^\vee)$.
If $\delta''= 0$, then $\Theta_{C/k}$ would have degree one, contradiction.
Hence $\delta''\neq 0$, and it follows $\shK^\vee\simeq\O_C(y_\infty)$.
\qed

\begin{remark}
\mylabel{no chance}
Let $\maxid\subset\O_Y$ be the maximal ideal for the point at infinity $y_\infty\in Y$.
The preceding result tells us that  any global vector field on $Y$,
the corresponding tangent vector in $\Theta_{Y/k}(y_\infty)\subset\Hom_k(\maxid/\maxid^2,k)$
is   tangent to the curve of singularities $C\subset Y$.
We conclude that our $\delta=PD_e$ is in some sense the best possible choice
when it comes to twisting in Section \ref{Twisted del Pezzo surfaces}.
\end{remark}

\begin{remark}
Suppose $\shF$ is a  locally free sheaf of rank $n$ on the projective plane $X=\PP^2$.
The restriction to any line $L\subset X$  splits into
a direct sum  of invertible sheafs $\shF_L\simeq \O_L(d_1)\oplus\ldots\oplus\O_L(d_n)$,
say with $d_1\leq\ldots\leq d_n$.
This sequence of integers is called the \emph{splitting type}
of $\shF$ along the line $L$.
The preceding results tell us:
The generic splitting type of $\shF=\nu^*(\Theta_{Y/k})$
is given by the sequence $(2,3)$. The generic splitting type degenerates
to the special splitting type $(1,4)$ on those lines $L\subset\PP^2$
whose image in $D\subset Y$ is   Cartier or equals $C$.
\end{remark}

\section{Twisted del Pezzo surfaces}
\mylabel{Twisted del Pezzo surfaces}

We keep the assumptions as in the preceding section,
such that $Y$ is a nonnormal del Pezzo surfaces, defined by an even glueing
polynomial $P=\alpha_2u^2+\alpha_0$.
Then we have a global vector field $\delta\in H^0(Y,\Theta_{Y/k})$
with $\delta\circ\delta=0$ given by the formula $\delta=PD_e$.
Such vector fields correspond to actions of the   group scheme $\alpha_2$, which is finite and
infinitesimal.

Recall that we have $\alpha_2=\Spec k[\epsilon]$ as a scheme. Its values on
$k$-algebras $R$ is the group
$\alpha_2(R)=\left\{f\in R\mid f^2=0\right\}$, with addition as group law.
The action $\alpha_2\times Y\ra Y$ is given by the formula
$$
\O_Y\lra k[\epsilon]\otimes_k\O_Y ,\quad s\longmapsto \delta(s)\epsilon\otimes s.
$$
A rational point $y\in Y$ is a fixed point for the $\alpha_2$-action if 
and only if $\delta(y)=0$ as a section of $\Theta_{Y/k}$, or equivalently
$\delta(\maxid_y)\subset\maxid_y$ as   derivation $\delta:\O_Y\ra\O_Y$.

As explained in Section \ref{Twisted forms},
any $\alpha_2$-torsor $T$ yields a twisted form $Y'=Y\wedge T$ of our nonnormal del Pezzo surface $Y$.
Note that the projections $Y\leftarrow Y\times T\ra Y'$ are universal homeomorphisms,
and we may identify points on $Y$ with points on $Y'$.
Any such twisted form $Y'$ is locally of complete intersection. Moreover, $\omega_{Y'}$
is antiample, and we  have $h^1(\O_{Y'})=h^1(\O_Y)=1$. Whence the twisted form
$Y'$ is another del Pezzo surface, possibly with less severe singularities than $Y$.

Any   $\alpha_2$-torsor is of the form $T=\Spec k(\sqrt{\lambda})$
for some scalar $\lambda\in k$, with action given by the derivation $\sqrt{\lambda}\mt 1$.
The torsor is nontrivial if and only if $\lambda\in k$ is not a square.
We now can formulate the  main result of this paper:

\begin{theorem}
\mylabel{twisted surface}
Let $k$ be a nonperfect field of characteristic two, $\lambda\in k$ be a nonsquare, and $T=\Spec k(\sqrt{\lambda})$
the corresponding $\alpha_2$-torsor.
Then the twisted form $Y'=Y\wedge T$ is a normal del Pezzo surface
with $h^1(\O_{Y'})=1$.
It has a unique singularity $y_\infty'\in Y'$, which  corresponds to
the point at infinity $y_\infty\in Y$.
\end{theorem}

\proof
First observe that $Y'$ is smooth outside the curve corresponding
to the reduced singular locus $C\subset Y$.
According to Proposition \ref{zeros field}, the point at infinity 
$y_\infty\in Y$ is a  fixed point on the singular locus.
The corresponding point on the twisted form $y'_\infty\in Y'$ then
must be a singularity, by Proposition \ref{singularities stay}.

It remains to see what the effect of twisting is on the affine
open subset $V\subset Y$ at the singular locus.
Recall that the open subset $V_{P^2}\subset V$ is given by the algebra
$
k[a,b,e,P^{-1}]/(P^4e^2+b^4a+a^7)
$,
and that $\delta=PD_e$, compare the proof for Theorem \ref{complete intersection}.
We first analyse the rational point $y\in Y$ corresponding
to the origin $a=b=e=0$. We have $\delta(a)=\delta(b)=0$ and $\delta(e)=P$,
whence the orbit $Gy\subset Y$ is given by the ideal $(a,b,e^2)$.
Clearly, this ideal is generated by $a,b$, which is a regular sequence.
Now Theorem \ref{twisted forms} tells us that $Y'$ is regular near the
point corresponding to $y\in Y$.

Next we treat the singular locus of $V_{P^2}$ outside the origin.
Consider the Cartier divisor $A\subset V_{P^2}$ supported by the singular locus given by the ring element
$c=v^2=(b^2+a^3)/P^2$. This element is invariant, and we have
$$
k[a,b,e,P^{-1}]/(P^4e^2+b^4a+a^7,b^2+a^3) = k[a,b,e,P^{-1}]/(b^2+a^3,e^2).
$$
Twisting this algebra, we obtain as twisted algebra $k(\sqrt{\lambda})[a,b,P^{-1}]/(b^2+a^3)$,
which defines a rational cuspidal curve over the  quadratic extension field $k(\sqrt{\lambda})$.
The latter is regular outside the origin.
Using Theorem \ref{twisted forms} again, we conclude that the twisted form
$Y'$ is regular on the open subset corresponding to $V_{P^2}\subset Y$.

Finally, we treat the other open subset $V_Q\subset V$, which is given by
$$
k[a,b,c,c',a^{-1}]/(b^2+a^3+cP^2,c'^2=c^3).
$$
Here our derivation takes the form $\delta=aD_{c'}$.
Again, consider the Cartier divisor $A\subset V_Q$ supported by the singular locus given
by $c$. We have
$$
k[a,b,c,c']/(b^2+a^3+cP^2,c'^2=c^3,c)=k[a,b,c']/(b^2+a^3,c'^2),
$$
and we may argue as above.
The upshot is that the twisted form $Y'$ is regular outside the point
at infinity $y_\infty'\in Y'$.
\qed

\medskip
It is not difficult to analyse the singularity:

\begin{theorem}
\mylabel{rational double}
The singularity $y_\infty'\in Y'$ is a rational double point of type $A_1$.
The minimal resolution of singularities $r:\tilde{Y}'\ra Y' $
is obtained by blowing up the reduced singular point.
The exceptional divisor $E=r^{-1}(y_\infty')$
is isomorphic to a regular quadric in $\PP^2_k$ that is a twisted form
of the double line. The regular surface $\tilde{Y}'$ is a weak del Pezzo surface
with $h^1(\O_{\tilde{Y}'})=1$.
\end{theorem}

\proof
As explained in the proof for Proposition \ref{complete intersection}, 
the completion $\O^\wedge_{Y,y_\infty}$ is the algebra $R=k[[x,y,z]]/(y^3+z^2)$, and furthermore $\delta=xD_z$.
It follows that the completion   $\O^\wedge_{Y',y_\infty'}$
is the subalgebra $R\subset k[\sqrt{\lambda},x,y,z]/(y^3+z^2)$ generated
by the invariants $x,y,z'$, where $z'=z+\sqrt{\lambda}x$, which has defining
relation $z'^2=y^3+\lambda x^2$.

Consider the blowing up $Z\ra\Spec(R)$ of the maximal ideal $(x,y,z')$.
It is covered by two charts: The $x$-chart
\begin{equation}
\label{x-chart}
x,y/x,z'/x \quad \text{modulo}\quad (z'/x)^2 = (y/x)^3x +\lambda,
\end{equation}
and the $z'$-chart 
\begin{equation}
\label{z'-chart}
x/z',y/z',z'  \quad \text{modulo}\quad z'^2 = (y/z')^3z' +\lambda(x/z')^2.
\end{equation}
The exceptional divisor $E\subset Z$ is given by setting $x$ and $z'$ to zero,
respectively. Whence $E$ is covered by $y/x,z'/x$ modulo $(z'/x)^2=\lambda$
and $x/z,y/z'$ modulo $(x/z')^2=1/\lambda$.
The exceptional divisor is evidently regular, and hence the blowing up $\tilde{Y}'$ is regular as well.
Furthermore, we easily infer that $E=r^{-1}(y_\infty')$ is isomorphic to a regular quadric in $\PP^2_k$ that becomes
a double line after adjoining $\sqrt{\lambda}$.
 
We infer  that $R^1r_*(\O_{\tilde{Y}'})=0$, so the singularity is rational. 
It also follows that the map $H^1(Y',\O_{Y'})\ra H^1(\tilde{Y}',\O_{\tilde{Y}'})$
is bijective. Finally, write the relative dualizing sheaf in the form $\omega_{\tilde{Y}'/Y'}=\O_{\tilde{Y}'}(nE)$
for some integer $n$. Using
\begin{equation}
\label{discrepancy}
-2 =\deg(\omega_E)=\omega_{\tilde{Y}'}(E)\cdot E = (n+1)E^2,
\end{equation}
we conclude $n=0$ and $E^2=-2$. In other words, the singularity is a rational double point
of type $A_1$. Moreover, we have $\omega_{\tilde{Y}'}=r^*(\omega_{Y'})$, and hence the
antidualizing sheaf for $\tilde{Y}'$ is nef and big. In other words, the regular surface
$\tilde{Y}'$ is a weak del Pezzo surface.
\qed

\medskip
According to Mumford's result \cite{Mumford 1961},   the only normal surface
singularities
over the complex numbers whose formal completion are factorial
are the rational double points of type $E_8$
(for arbitrary algebraically closed ground fields, see \cite{Lipman 1969}, \S 25).
The situation is more complicated over nonclosed ground fields.  

\begin{corollary}
\mylabel{locally factorial}
The complete local rings $\O^\wedge_{Y',y'}$ of our twisted del Pezzo surface $Y'$ are
factorial.
\end{corollary}

\proof
Let $D\subset Y'$ be a Weil divisor, and $\tilde{D}\subset \tilde{Y}'$ be
its strict transform.
The exceptional divisor $E\subset\tilde{Y}'$ carries no invertible sheaf of
degree one. Rather, it is a cyclic group generated by the invertible sheaf $\O_E(E)$,
which has degree two.
Write  $\tilde{D}\cdot E=2n$ for some integer $n$. Then $(\tilde{D}+nE)\cdot E=0$.
This implies that the invertible sheaf $\shL=\O_{\tilde{Y}'}(\tilde{D}+nE)$
is trivial on the formal
 completion along $E$, because  $H^1(\tilde{Y}',\O_{mE})=0$ for all integers $m\geq 0$.
It follows that the coherent $\O_{Y'}$-module $r_*(\shL)$  is invertible.
Therefore, the Weil divisor $D\subset Y'$ must be Cartier.
The same argument applies for formal Weil divisors on   $\Spec(\O^\wedge_{Y',y'})$.
\qed

\section{Fano-Mori contractions}
\mylabel{Fano--Mori contractions of fiber type}

We now use the results of the preceding section on del Pezzo surfaces over nonperfect
ground fields to construct some interesting Fano-Mori contractions of fiber type over
algebraically closed fields.
Let us  now work, for simplicity, over an algebraically closed ground field $k$ of characteristic two.

Choose an abelian variety $A'$  with $a$-number $a(A)\geq 1$.
This mean that there exists at least one embedding $\alpha_2\subset A'$, and in turn an
$\alpha_2$-action on  the abelian variety via translations.
In dimension one, for example, we could choose a supersingular elliptic curve
with Weierstrass equation of the form $y^2+y= x^3+a_4x+a_6$, with action given by the derivation
$x\mt 1$, $y\mt x^2+a_4$. Note that in characteristic two, all supersingular elliptic
curves are isomorphic.
The quotient $A=A'/\alpha_2$ is again an abelian variety, and
the quotient map $A'\ra A$ is a purely inseparable isogeny of degree two.

We now fix once and for all an embedding $\alpha_2\subset A'$ and consider the
corresponding $\alpha_2$-action on $A'$ via translations.
Let $Y$ be the nonnormal del Pezzo surface constructed in the preceding sections.
We assume that the glueing polynomial $P$ is even, such that we have
the global vector field $\delta=PD_e$ corresponding to an $\alpha_2$-action on $Y$.
The product $Z'=Y\times A'$ carries the diagonal action, and we may take the
quotient $Z=\alpha_2\backslash Z'$.
The projection $f':Z'\ra A'$ induces a projection $f:Z\ra A$.
To understand its fibers, consider the function fields $K'=k(A')$
and $K=k(A)$. Then $K\subset K'$ is a purely inseparable quadratic field extension,
and hence of   the form 
$K'=K(\sqrt{\lambda})$ for some nonsquare $\lambda\in K$.
We may view $T=\Spec K'$ as an $\alpha_2$-torsor over $K$.

\begin{proposition}
The generic fiber $Z_\eta$ of the projection $f:Z\ra A$
is the twisted form $Y_K\wedge T$, which is a normal del Pezzo surface.
For all closed points $\sigma\in A$, the fiber
$Z_\sigma$ is isomorphic to the nonnormal del Pezzo surface $Y$.
\end{proposition}

\proof
Taking quotients by free group actions commutes with arbitrary base change.
Given a point $\sigma\in A$ with residue field $\kappa=\kappa(\sigma)$, and
$T\subset A'$ be its preimage.
Making base change with respect to $T\ra A$, we see that
the  fiber $Z_\sigma$ is the quotient of $Y_{\kappa}\times_{\Spec(\kappa)} T$ by the diagonal action, 
so that $Z_\sigma= Y_\kappa\wedge T$.
If $\sigma$ is a closed point, the torsor $T$ is trivial, and hence $Z_\sigma=Y$.
If $\sigma$ is the generic point, then the torsor $T$ is nontrivial.
According to Theorem \ref{twisted surface}, the twisted form is then normal.
\qed

\medskip
The point at infinity $y_\infty\in Y$ is invariant under
the $\alpha_2$-action.
Whence it defines a section $s:A\ra Z$, whose image is
the quotient of $\left\{y_\infty\right\}\times A'$ by the diagonal action.

\begin{proposition}
\mylabel{normal complete intersection}
The scheme $Z$ is normal  and locally of complete intersection.
The reduced singular locus of $Z$ equals the image of the section 
$s(A)\subset Z$.
\end{proposition}

\proof
The  morphism $f:Z\ra A$ is flat, because the composition
$Z'\ra A$ is flat and the quotient $Z'\ra Z$ is faithfully flat.
The base $A$ and all fibers $Z_a$ are locally of complete intersection,
whence $Z$ is locally of complete intersection.

According to Theorem \ref{twisted surface}, the singular locus of  the generic fiber
ist $s(A)_\eta$.
It follows that  $s(A)\subset\Sing(Z)$.
The translation action of $A'$ on $Z'=Y\times A'$ via the second factor
commutes with the diagonal $\alpha_2$-action, whence induces
an action of $A$ on the quotient $Z$.
For any closed point $z\in Z$, the induced map $f:Az\ra A$ is surjective.
Using that the singular locus is invariant under this action,
we infer that $\Sing(Z)\subset s(A)$.
In particular, $Z$ is regular in codimension one. It follows
that $Z$ is normal.
\qed

\begin{proposition}
\mylabel{canonical singularities}
The blowing up $r:\tilde{Z}\ra Z$ with center the reduced subscheme
$s(A)\subset Z$ is a resolution of singularities.
The singularities of $Z$ are canonical.
\end{proposition}

\proof
According to \ref{rational double}, the generic fiber
$\tilde{Z}_\eta$ is regular. Moreover, the $A$-action on $Z$
leaves the center of the blowing up $s(A)\subset Z$ invariant,
and hence the action extends to the relative homogeneous spectrum $\tilde{Z}=\shProj(\oplus(\shI^n/\shI^{n+1}))$,
where $\shI\subset\O_Z$ denotes the ideal of the center.
We now may argue as in the preceding proof and infer that $\tilde{Z}$ must be regular.

Let $E\subset\tilde{Z}=r^{-1}(s(A))$ be the exceptional divisor.
It must be flat over the base of the projection $\tilde{Z}\ra A$ because
it carries an $A$-action.
The relative dualizing sheaf $\omega_{\tilde{Z}/Z}$ is of the form
$\O_{\tilde{Z}}(nE)$ for some integer $n$, which is called the \emph{discrepancy}
for the resolution of singularities. The singularities on the threefold $Z$ are called \emph{canonical}
if $n\geq 0$. According to (\ref{discrepancy}), we have $n=0$, and hence $Z$ is canonical.
\qed

\medskip
Recall that a   morphism of proper normal scheme $f:V\ra W$ is called a \emph{Fano-Mori contraction}
if $\O_W\ra f_*(\O_V)$ is bijective, the total space $V$ is $\QQ$-Gorenstein,
and $\omega_V^\vee$ is $f$-ample.

\begin{proposition}
\mylabel{singular locus}
The morphism $f:Z\ra A$ is a Fano-Mori contraction. The
$\O_A$-module $R^1f_*(\O_Z)$ is invertible and commutes with base change.
\end{proposition}

\proof
For all closed points $a\in A$, we have 
$Z_a=Y$, and hence $h^0(\O_{Z_a})=1$ and $h^2(\O_{Z_a})=0$.
If follows that the canonical map $\O_A\ra f_*(\O_Z)$ is bijective,
and that the coherent $\O_A$-module $R^1f_*(\O_Z)$ is locally free, of
rank $h^1(\O_{Z_a})=1$.

By Proposition \ref{normal complete intersection}, the scheme $Z$ is Gorenstein.
Let $C\subset Z$ be an integral curve, and $C'\subset Z'$ be its
preimage.
We have $\omega_Z\otimes\O_{Z'}=\pr_1^*(\omega_Y)$.
If follows that $C\cdot \omega_Z\leq 0$, with equality if and only if
the induced map
$f:C\ra A$ is finite. Summing up, $f:Z\ra A$ is a Fano-Mori contraction.
\qed

\section{Maps to projective spaces}
\mylabel{Maps to projective spaces}

We now return to our nonnormal del Pezzo surface  $Y$.
In this final section we study maps to  projective spaces, which are defined in terms of
semiample invertible sheaves.  The upshot will be that it is not possible
to define $Y$ in  a simple way as a hypersurface in   projective space, or a finite covering
of  projective space. Obviously, the same then holds for   twisted forms $Y'=Y\wedge T$.

The various geometric properties of invertible sheaves on $Y$
can be nicely expressed in terms of degrees.

\begin{theorem}
\mylabel{semiample degree}
An invertible $\O_Y$-module $\shL\neq\O_Y$ of degree $d=\deg(\shL)$ is$:$
\renewcommand{\labelenumi}{(\roman{enumi})}
\begin{enumerate}
\item
 semiample if and only if $d\geq 0$\! $;$
\item
 ample if and only if $d\geq 1$\! $;$
\item
globally generated if and only if $d\geq 2$\! $;$
\item
very ample if and only if $d\geq 3$\! $;$
\end{enumerate}
\end{theorem}

\proof
For this we may assume that the ground field $k$ is algebraically closed.
We already proved assertion (ii) in Proposition \ref{ample degree}.
Concerning (i), recall that semiampleness means that some tensor power is globally generated.
Suppose $\shL$ is semiample. Then the restriction $\shL_C$ is semiample as well,
and hence $d\geq 0$. Conversely, suppose $d\geq 0$. If $d>1$, then $\shL$ is ample,
and if $d=0$, then the sheaf  $\shL^{\otimes 2}$ is trivial by the exact sequence (\ref{picard extension}).
In both cases $\shL$ is semiample.

Next, we prove (iii). Suppose that $\shL$ is globally generated.
Then the restriction $\shL_C$ to the cuspidal curve of arithmetic genus $p_a=1$ 
is globally generated as well, and
this implies that $d\geq 2$ by Lemma \ref{hyperelliptic curve} below.
Conversely, suppose the degree is $d\geq 2$.
Decompose $\shL=\shL_1\otimes\cdots\otimes\shL_d$ into a tensor product
of invertible sheaves of degree one.
According to Proposition \ref{euler characteristic}, 
we have $h^0(\shL_i)=1$, whence there are unique Cartier divisors
$C_i\subset Y$ with $\shL_i=\O_Y(C_i)$.
It follows that the base locus of $\shL$ is contained in $\bigcup_{i=1}^d C_i$.
The exact sequence
$0\ra\shL(-C_i)\ra\shL\ra\shL_{C_i}\ra 0$
yields an exact sequence
$$
H^0(Y,\shL)\lra H^0(C_i,\shL_{C_i})\lra H^1(Y,\shL(-C_i)).
$$
The term on the right vanishes by Proposition \ref{euler characteristic},
because $\shL(-C_i)$ has degree $d-1\geq 1$.
To finish the argument, it suffices to check that $\shL_{C_i}$
is globally generated. According to Proposition \ref{cartier weil}, the
$C_i$ are isomorphic to the rational cuspidal curve with $p_a=1$.
By Lemma \ref{hyperelliptic curve} below, $\shL_{C_i}$ is globally generated.

It remains to prove (iv), which is the most interesting part.
Suppose first that $\shL$ is very ample. Then the restriction $\shL_C$
is very ample as well, and this implies $d\geq 3$ by Lemma 
\ref{hyperelliptic curve}
below.
Conversely, suppose $d\geq 3$. Let $A\subset Y$ be an Artin subscheme
of length two.
We have to show that $H^0(Y,\shL)\ra H^0(A,\shL_A)$ is surjective.
The idea is to use Cartier divisors of degree two.

Let $\shN$ be an invertible $\O_Y$-module of degree two.
We already saw that $\shN$ is globally generated, whence
defines a morphism $r_\shN:Y\ra\PP^2$.
The image $r_\shN(A)\subset\PP^2$ is an Artin scheme of length $\leq2$,
and hence $\shN$ has a nonzero global section whose zero scheme
$D\subset Y$ contains $A$.
The exact sequence $0\ra\shL(-D)\ra\shL\ra\shL_D\ra 0$ yields an exact sequence
$$
H^0(Y,\shL)\lra H^0(D,\shL_D)\lra H^1(Y,\shL(-D)).
$$
The term on the right vanishes by Proposition \ref{euler characteristic}, because $\shL(-D)$ has degree $\geq 1$.
Hence it suffices to show that $H^0(D,\shL_D)\ra H^0(A,\shL_A)$ is surjective.

Now suppose for a moment that $\omega_Y\otimes\shN^{\otimes 2}\neq \shL$
and that $r(A)\subset \PP^2$ has length one.
Using the latter, we see that $\shN$ has another nonzero section
whose zero scheme $D'\subset Y$ having no irreducible component in common with $D$
and containing $A$. Set $A'=D\cap D'$. Then $A\subset A'$,
and the inclusion $A'\subset D$ is Cartier.
The exact sequence $0\ra\shL_D(-A')\ra\shL_D\ra\shL_{A'}\ra0$
yields an exact sequence
$$
H^0(D,\shL_D)\lra H^0(A',\shL_{A'})\lra H^1(D,\shL_D(-A')).
$$
The term on the right sits  inside the exact sequence
$$
H^1(Y,\shL\otimes\shN^\vee)\lra H^1(D,\shL_D(-A'))\lra 
H^2(Y,\shL\otimes\shN^{\otimes -2}).
$$
In this sequence, the term on the left vanishes by Proposition \ref{euler characteristic}, since we have $\deg(\shL\otimes\shN^\vee)<0$. The term on the right is Serre dual
to $H^0(Y,\shM)$,
where $\shM=\omega_Y\otimes\shL^{\otimes -1}\otimes\shN^{\otimes 2}$.
This sheaf has degree $3-d$, and hence $H^0(Y,\shM)$ vanishes for $d>3$. In the boundary case
$d=3$ we also have $H^0(Y,\shM)=0$, because we are presently assuming
that $\omega_Y\otimes\shN^{\otimes 2}\neq \shL$.
Combining these observations, we see that the restriction map
$H^0(Y,\shL)\ra H^0(A,\shL_A)$ is surjective.

To complete the proof, we now may assume that for all
invertible sheaves $\shN$ of degree two with 
$\omega_Y\otimes\shN^{\otimes 2}\neq \shL$ the image $r_\shN(A)\subset\PP^2$
has length two.
Our goal now is to find a global section of $\shL$
whose zero scheme   intersects $A$ but does not contain  $A$. This implies that
$H^0(Y,\shL)\ra H^0(A,\shL_A)$ is surjective, because
we already know that $\shL$ is globally generated.
By our assumption, for any $\shN$ of degree two with
$\omega_Y\otimes\shN^{\otimes 2}\neq \shL$, we  find a global section of $\shN$
whose zero scheme $D\subset Y$ intersects $A$ but does not contain  $A$.

We now have to distinguish the cases that $n=\deg(\shL)$
is even or odd. I only go through the case that $n=2m+1$ is odd, the even
case being similar. Choose a Cartier divisor $D'\subset Y$ of degree two disjoint
from $A$. The  equation of invertible sheaves $\shL\simeq \shN((m-1)D'+E)$
defines an invertible sheaf $\O_Y(E)$ of degree one.
Since $H^0(Y,\O_Y(E))=1$, the effective Cartier divisor $E\subset Y$
is also unique.
Now note that  $E\subset Y$ is the image of a line $L\subset\PP^2$,
and that the restriction map $H^0(Y,\O_Y(E))\ra H^0(C,\O_C(E))$
is bijective, according to Proposition \ref{restriction bijective}.
From this we infer that there is at most one invertible sheaf of degree
one $\O_Y(E)$ with $A\subset E$.
So tensoring $\shN$ with some general numerically trivial invertible sheaf,
we may assume that $A\not\subset E$.
If $A$ is disjoint from $E$, then 
$D + (m-1)D' + E$ is the desired Cartier divisor representing $\shL$
that intersect but does not contain $A$.
If $A\cap E$ is nonempty, we simply replace $D$ by a linearly equivalent
Cartier divisor disjoint form $A$, and conclude as above.
\qed

\medskip
Suppose that the invertible $\O_Y$-module $\shL$ is globally generated.
In other words, its degree is $d\geq 2$.
Set $n=d(d+1)/2-1$, and let $r_\shL:Y\ra\PP^n$ be the resulting morphism
defined by $\shL$.

\begin{corollary}
\mylabel{minimal maps}
If $d=2$, then the morphism $r_\shL:Y\ra\PP^2$ is flat, surjective, of degree four, and
all fibers are Artin schemes  of complete intersection.
There is no surjection to the projective plane of smaller degree.
If $d=3$, then $r_\shL$ is a closed embedding $Y\subset\PP^5$.
There is no closed embedding into any projective space of smaller dimension.
\end{corollary}

\proof
Suppose $d=2$. The morphism $r_\shL:Y\ra\PP^2$ is flat because $Y$ is Cohen-Macaulay and
$\PP^2$ is regular (\cite{Serre 1965}, page IV-37, Proposition 22).
The other statements follow immediately from Theorem \ref{semiample degree}.
\qed

\medskip

In the course of the proof for Theorem \ref{semiample degree}, we used the following facts.

\begin{lemma}
\mylabel{hyperelliptic curve}
Let $C$ be the rational cuspidal curve of arithmetic genus $p_a=1$, and $\shL\neq\O_C$ be an invertible
$\O_C$-module of degree $d$. Then $\shL$ is globally generated if and only if $d\geq 2$,
and very ample if and only if $d\geq 3$.
\end{lemma}

\proof
Of course, we may assume that the ground field $k$ is algebraically closed.
The arguments are similar to the case of elliptic curves. The problem, however,
is that some Weil divisors on $C$ are not Cartier.

Let us first prove that the numerical conditions are necessary.
Suppose that $\shL$ is globally generated, so $d\geq 0$.
The case $d=0$ is  impossible, because $\shL$ is nontrivial by assumption.
Hence $d\geq 1$, and the usual argument gives $h^0(C,\shL)=d$.
The invertible sheaf $\shL$ is ample,  whence the morphism $r_{\shL}:C\ra\PP^{d-1}$
is finite, and therefore $d\geq 2$.
If, furthermore, $\shL$ is very ample, we must have $d\geq 3$.

The converse is more interesting.
Suppose $d\geq 2$, and let $y\in C$ be any closed point, which
is a Weil divisor of length one. I claim that there is a 
Cartier divisor $D\subset C$ of degree at most two that contains $y$
and has $\shL\not\simeq\O_C(D)$.
Suppose this for the moment. 
The short exact sequence $0\ra\shL(-D)\ra\shL\ra\shL_{D}\ra 0$ yields an exact
sequence
\begin{equation}
\label{curve sequence}
H^0(C,\shL)\lra H^0(D,\shL_{D})\lra H^1(C,\shL(-D)).
\end{equation}
The term on the right is Serre dual to $H^0(C,\shL^\vee(D))$.
The invertible sheaf $\shL^\vee(D)$ has degree $n\leq 2-d\leq 0$,
and in the case $n=0$ is nontrivial. Whence it has no global section,
and it follows that $\shL$ has a section that does not vanish at $y$.

Let us now verify the claim. There is nothing to prove if $y\in C$
is contained in the regular locus. So let us assume that it is the
singular point, and write $\O_{C,y}^\wedge=k[[u^2,u^3]]$.
For any scalar $\lambda\in k$, the element $u^2+\lambda u^3$ defines
a Cartier divisor $D_\lambda\subset C$ of length two with support $y$.

Finally, suppose that $d\geq 3$.
Let $A\subset C$ be a closed subset of length two.
We have to see that the restriction map $H^0(C,\shL)\ra H^0(A,\shL_A)$ is surjective.
It follows from the above that there is a Cartier divisor $D\subset A$ of length four
containing $A$, and with $\shL\not\simeq\O(D)$.
In the case $d\geq 4$, one proceeds easily as above to see that $H^0(C,\shL)\ra H^0(A,\shL_A)$ is surjective.
For $d=3$ we argue as follows: We have $h^0(C,\shL)=3$, and we already know
that $\shL$ is globally generated, hence there is a finite
morphism $r_\shL:C\ra\PP^2$, which does not factor over
a line $\PP^1\subset\PP^2$. Whence $r_\shL$ is birational onto its
image $r_\shL(C)$, which must be a cubic. Any cubic has arithmetic genus $p_a=1$.
Since $C$ also has arithmetic genus $p_a=1$, the birational morphism $r_\shL$
must be an isomorphism.
\qed



\begin{thebibliography}{99}


\bibitem{Artin 1970}
\textsc{M.~Artin},
Algebraization of formal moduli II: Existence of modifications,
Ann.\ of Math.\ (2) 91 (1970), 88--135.

\bibitem{Avramov 1975}
\textsc{L.\ Avramov},
Flat morphisms of complete intersections,
Dokl.\ Akad.\ Nauk SSSR  225  (1975),  11--14.

\bibitem{Berthelot; Bloch; Esnault 2005}
\textsc{P.\ Berthelot, S.\ Bloch and H.\ Esnault},
On Witt vector cohomology for singular varieties, 
Preprint, math.AG/0510349.

\bibitem{Bombieri; Mumford 1977}
\textsc{E.~Bombieri and D.~Mumford},
Enriques' classification of surfaces in char.\ $p$.  II,
Complex analysis and algebraic geometry, pp.\ 23--42,
Iwanami Shoten, Tokyo, 1977.

\bibitem{Bombieri; Mumford 1976}
\textsc{E.~Bombieri and D.~Mumford},
Enriques' classification of surfaces in char $p$.  III,
Invent.\ Math.\ 35  (1976), 197--232.

\bibitem{Bayer; Eisenbud 1995}
\textsc{D.\ Bayer and D.\ Eisenbud},
Ribbons and their canonical embeddings,
Trans.\ Am.\ Math.\ Soc.\ 347 (1995), 719--756.

\bibitem{Demazure 1976}
\textsc{M. Demazure},
Surfaces de Del Pezzo,
Seminaire sur les singularites des surfaces, pp.\ 21--70,
Lecture Notes in Math.\  777,
Springer, Berlin, 1980.

\bibitem{Esnault 2003}
\textsc{H.\ Esnault},
Varieties over a finite field with trivial Chow group of 0-cycles have a rational point,  
Invent.\ Math.\  151  (2003),   187--191.

\bibitem{Giraud 1971}
\textsc{J.~Giraud},
Cohomologie non ab\'elienne, 
Grundlehren Math.\ Wiss.\ 179, Springer, Berlin, 1971.

\bibitem{EGA II}
\textsc{A.\ Grothendieck},
\'El\'ements de g\'eom\'etrie alg\'ebrique II,
\'Etude globale \'el\'ementaire de quelques classes de morphismes,
Inst.\ Hautes \'Etud.\ Sci.\ Publ.\ Math.\ 8 (1961).

\bibitem{EGA IIIa}
\textsc{A.\ Grothendieck},
\'El\'ements de g\'eom\'etrie alg\'ebrique III,
\'Etude cohomologique des faiscaux coh\'erent,
Inst.\ Hautes \'Etud.\ Sci.\ Publ.\ Math.\ 11 (1961).

\bibitem{EGA IVa}
\textsc{A.\ Grothendieck},
\'El\'ements de g\'eom\'etrie alg\'ebrique IV, \'Etude locale des
sch\'emas et des morphismes de sch\'emas,
Inst.\ Hautes \'Etud.\ Sci.\ Publ.\ Math.\ 20 (1964).

\bibitem{EGA IVc}
\textsc{A.\ Grothendieck},
\'El\'ements de g\'eom\'etrie alg\'ebrique IV, \'Etude locale des
sch\'emas et des morphismes de sch\'emas,
Inst.\ Hautes \'Etud.\ Sci.\ Publ.\ Math.\  28 (1966).

\bibitem{EGA IVd}
\textsc{A.\ Grothendieck},
\'El\'ements de g\'eom\'etrie alg\'ebrique IV, \'Etude locale des
sch\'emas et des morphismes de sch\'emas,
Inst.\ Hautes \'Etud.\ Sci.\ Publ.\ Math.\   32 (1967).

\bibitem{SGA 1}
\textsc{A.~Grothendieck et al.},
Rev\^etements \'etales et groupe fondamental,
Lecture Notes in Math.\  224,
Springer, Berlin, 1971.

\bibitem{SGA 3a}
\textsc{A.~Grothendieck et al.},
Schemas en groupes I,
Lecture Notes in Math.\  151,
Springer, Berlin, 1970.

\bibitem{Hidaka; Watanabe 1981}
\textsc{F.~Hidaka and K.~Watanabe},
Normal Gorenstein surfaces with ample anti-canonical divisor,
Tokyo J.\ Math.\ 4  (1981), 319--330.

\bibitem{Kawamata 1982}
\textsc{Y.\ Kawamata},
A generalization of Kodaira--Ramanujam's vanishing theorem,  
Math.\ Ann.\  261  (1982),   43--46.

\bibitem{Keel; Mori 1997}
\textsc{S.\ Keel and S.\ Mori},
Quotients by groupoids, 
Ann.\ of Math.\ (2) 145 (1997),  193--213.

\bibitem{Kodaira  1953}
\textsc{K.\ Kodaira},
On a differential-geometric method in the theory of analytic stacks,
Proc.\ Nat.\ Acad.\ Sci.\ U.S.A.\ 39 (1953), 1268--1273.

\bibitem{Kollar 1991}
\textsc{J.\ Koll\'ar},
Extremal rays on smooth threefolds,
Ann.\ Sci.\ \'Ecole Norm.\ Sup.\ (4) 24  (1991), 339--361.

\bibitem{Lauritzen 1996}
\textsc{N.\ Lauritzen},
Embeddings of homogeneous spaces in prime characteristics,    
Amer.\ J.\ Math.\  118  (1996), 377--387.

\bibitem{Lauritzen; Rao 1997}
\textsc{N.\ Lauritzen and A.\ Rao},
Elementary counterexamples to Kodaira vanishing in prime characteristic,  
Proc.\ Indian Acad.\ Sci.\ Math.\ Sci.\  107  (1997),  21--25.

\bibitem{Lipman 1969}
\textsc{J.\ Lipman:}
Rational singularities, with applications to algebraic surfaces and unique factorization, 
Inst.\ Hautes \'Etud.\ Sci.\ Publ.\ Math.\  36 (1969), 195--279.

\bibitem{Megyesi 1998}
\textsc{G.\ Megyesi},
Fano threefolds in positive characteristic,
J.\ Algebraic Geom.\  7  (1998), 207--218.

\bibitem{Mori; Saito 2003}
\textsc{S.\ Mori and N.\ Saito},
Fano threefolds with wild conic bundle structures,
Proc.\ Japan Acad. Ser.\ A Math.\ Sci.\  79  (2003), 111--114.

\bibitem{Mumford 1961}
\textsc{D.~Mumford},
The topology of normal singularities of an algebraic surface and a criterion for simplicity,
Inst.\ Hautes \'Etud.\ Sci.\ Publ.\ Math.\ 9 (1961), 5--22.

\bibitem{Mumford 1966}
\textsc{D.~Mumford},
Lectures on curves on an algebraic surface,
Ann.\ of Math.\ Stud.\ 59, Princeton University Press, Princeton, N.J., 1966.

\bibitem{Raynaud 1978}
\textsc{M.\ Raynaud},
Contre-exemple au ``vanishing theorem" en caract\'eristique $p>0$, 
C.\ P.\ Ramanujam---a tribute, pp. 273--278,
Tata Inst.\ Fund.\ Res.\ Studies in Math.\ 8.
Springer, Berlin, 1978. 
 
\bibitem{Reid 1994}
\textsc{M.\ Reid},
Nonnormal del Pezzo surfaces,
Publ.\ Res.\ Inst.\ Math.\ Sci.\  30  (1994), 695--727.

\bibitem{Saito 2003}
\textsc{N.\ Saito},
Fano threefolds with Picard number 2 in positive characteristic,
Kodai Math.\ J.\  26  (2003), 147--166.

\bibitem{Schroeer 2001}
\textsc{S.~Schr\"oer},
Normal del Pezzo surfaces containing a nonrational singularity,
Manuscr.\ Math.\ 104 (2001), 257--274.

\bibitem{Schroeer; Siebert 2002}
\textsc{S.\ Schr\"oer and B.\ Siebert},
Irreducible degenerations of primary Kodaira surfaces,
Complex Geometry (G\"ottingen 2000),
pp.\ 193--222,
Springer, Berlin, 2002.

\bibitem{Schroeer 2005}
\textsc{S.\ Schr\"oer},
Kummer surfaces  for the selfproduct of the cuspidal rational curve,
math.AG/0504023, to appear in J. Algebraic Geom.

\bibitem{Serre 1965}
\textsc{J.-P.\ Serre},
Alg\`ebre locale. Multiplicit\'es, 
Lecture Notes  in Math.\ 11,
Springer, Berlin, 1965.

\bibitem{Serre 1972}
\textsc{J.-P.\ Serre},
Cohomologie galoisienne,
Fifth edition. Lect.\ Notes  Math.\ 5, 
Springer, Berlin, 1994.

\bibitem{Serre 1975}
\textsc{J.-P.~Serre},
Groupes alg\'ebriques et corps de classes,
Deuxi\`eme \'edition, Actualit\'es Scientifiques et Industrielles 1264,
Hermann, Paris, 1975.

\bibitem{Shepherd-Barron 1997}
\textsc{N.\ Shepherd-Barron},
Fano threefolds in positive characteristic,
Compositio Math.\  105  (1997),  237--265.

\bibitem{Viehweg 1982}
\textsc{E.\ Viehweg},
Vanishing theorems,
J.\ Reine Angew.\ Math.\ 335 (1982), 1--8.
\end{thebibliography}
\end{document}